\documentclass[fleqn,12pt]{mypaper}

\usepackage{times}
\usepackage[top=0.8in, bottom=0.8in, left=0.75in, right=0.75in]{geometry}
\usepackage{color}
\usepackage{algorithmic}
\usepackage[colorlinks=true,linkcolor=cyan]{hyperref}
\usepackage{upgreek}
\usepackage{float}
\usepackage{algorithm}

\def\equationautorefname~#1\null{%
  (#1)\null%
}

\let\oldnl\nl
\newcommand{\nonl}{\renewcommand{\nl}{\let\nl\oldnl}}

\usepackage[small]{titlesec}
\titleformat{\section}
  {\normalfont\fontsize{13}{15}\bfseries}{\thesection}{1em}{}
\titleformat{\subsection}
  {\normalfont\fontsize{12}{14}\bfseries}{\thesubsection}{1em}{}
\titleformat{\subsubsection}
  {\normalfont\fontsize{12}{14}\itshape}{\thesubsubsection}{1em}{}

\usepackage{soul}
\usepackage{hyperref}
\usepackage{graphicx}
\usepackage{bbm}  
\usepackage{bm}  
\usepackage{amsmath,amsthm,amssymb,amscd}
\usepackage{caption,subcaption}
\usepackage{booktabs,multirow}
\usepackage{setspace}
\usepackage[inline,final]{showlabels}
\usepackage{dashbox}
\usepackage{microtype}
\usepackage[shortcuts]{extdash}
\usepackage{paralist}
\usepackage{srcltx}

\newtheorem{proposition}{Proposition}[section]
\newtheorem{theorem}{Theorem}[section]

\newtheorem{observe}{Observation}[section]
\newtheorem{remark1}[observe]{Remark}

\usepackage{bm}

\begin{document}





\begin{center}
    \begin{minipage}[t]{6.0in}
    
 \begin{center}
   \begin{minipage}[t]{4.4in}
     \begin{center}

 \textbf{A high-order deferred correction method for the solution of free boundary
problems using penalty iteration, with an application to American option pricing} \\

  \vspace{ 0.50in}

 Dawei Wang$\mbox{}^{\dagger}$,
 Kirill Serkh$\mbox{}^{\ddagger\, \diamond}$, and
 Christina Christara$\mbox{}^{\dagger\, \circ}$  \\
 {\tt\small dwang@cs.toronto.edu, kserkh@math.toronto.edu, ccc@cs.toronto.edu}\\
\today

     \end{center}
   \vspace{ -100.0in}
   \end{minipage}
 \end{center}

 \vspace{ 2.00in}

This paper presents a high-order deferred correction algorithm combined with penalty iteration for solving free and moving boundary problems, using a fourth-order finite difference method. Typically, when free boundary problems are solved on a fixed computational grid, the order of the solution is low due to the discontinuity in the solution at the free boundary, even if a high-order method is used. Using a detailed error analysis, we observe that the order of convergence of the solution can be increased to fourth-order by solving successively corrected finite difference systems, where the corrections are derived from the previously computed lower order solutions.
The corrections are applied solely to the right-hand side, and leave the finite difference matrix unchanged.
The penalty iterations converge quickly
given a good initial guess. We demonstrate the accuracy and efficiency of our algorithm using several examples. Numerical results show that our algorithm gives fourth-order convergence for both the solution and the free boundary location.
We also test our algorithm on the challenging American put option pricing problem. Our algorithm gives the expected high-order convergence.

 \thispagestyle{empty}

   \vspace{ -100.0in}

   \end{minipage}
 \end{center}

 \vspace{ 2.60in}
 \vspace{ 0.50in}






 \vspace{ 2.00in}

 \vfill

\noindent
$\mbox{}^{\diamond}$  This author's work was supported in part by the NSERC
Discovery Grants RGPIN-2020-06022 and DGECR-2020-00356.
 \\
$\mbox{}^{\circ}$   This author's work was supported in part by the NSERC
Discovery Grant RGPIN-2021-03502.
 \\

 \vspace{2mm}

 \noindent
 $\mbox{}^{\dagger}$ Department of Computer Science, University of Toronto,
 Toronto, ON M5S 2E4\\
 \noindent
 $\mbox{}^{\ddagger}$ Departments of Mathematics and Computer Science, University of Toronto,
 Toronto, ON M5S 2E4 \\
 \vspace{2mm}


 \vfill
 \eject

\section{Introduction}\label{sec:intro}

Free boundary problems, in which both the solution to a partial differential
equation (PDE) and the domain on which it is defined are unknowns to be solved
for, arise in numerous applications of practical importance.  Many free
boundary problems, both with boundaries that move over time (called moving
boundary problems) and boundaries that are invariant with time, can be
reformulated as linear complementarity problems (LCPs) (see, for
example,~\S8.5 of~\cite{crank1984free}), a few well-known examples being the elliptic obstacle problem (see, for
example,~\cite{rodrigues1987obstacle}) and the American option pricing
problem (see, for example,~\cite{tavella2000pricing}).

Two important categories for solving free boundary problems are the front tracking methods and the fixed domain methods.
Front
tracking methods directly compute an approximation to the free boundary,
either at each time step in time-dependent problems, or iteratively in
%
%
time-independent problems (as in, for example,~\cite{sharma2019shape}).
While the free boundary can be tracked parametrically or as an indicator
function of some set, the most common
approaches are the level set method (see, for
example,~\cite{sethian2003level} for a survey), in which the free boundary
is represented as the zero level-set of a function which obeys an evolution
equation, and the phase-field method, in which the free boundary is
approximated by a finite-width region where a
phase-field function smoothly changes sign across the region (see, for
example,~\cite{anderson1998diffuse} for a survey). The front-fixing method,
in which the free boundary PDE is transformed into a nonlinear PDE with a fixed boundary,
is also considered to be a front tracking
method (see, for example,~\cite{wu1997front}).  Whatever the particular
method used, front tracking requires the construction of a separate
algorithm for approximating the front, derived from the underlying equations
and constraints of the free boundary problem.

In contrast, fixed domain methods reformulate the problem over the whole of
a fixed domain, and solve the new equations in such a way that the position
of the free boundary is returned simultaneously with the solution to the
PDE, and appears a posteriori as part of the solution process.  Such methods
have a reputation for being robust and relatively straightforward to
implement.
The most widely-used fixed-domain method is the penalty method, which
incorporates the inequality constraint of the LCP into the PDE by adding a nonlinear
penalty term (see, for example,~Ch. 1, \S8 of~\cite{friedman}). The
resulting penalized equation can be solved by successive over-relaxation
(SOR), which can be fairly expensive (see, for
example,~\cite{dewynne1993some}); this process can be accelerated by the
multigrid method (see, for
example,~\cite{clarke1999multigrid}).
%
%
Remarkably, under certain conditions, when
Newton's method is applied to the penalized PDE, the solution converges
monotonically and exhibits the rapid convergence characteristic of Newton's
method.
However, while the discretized
equation is solved to high accuracy, the approximation of the discrete
solution to the true solution of the LCP is of low order, due to
the disagreement between the free boundary and the fixed computational
grid.

The problem of reconciling a nonconforming boundary with a fixed
computational grid has been studied extensively, particularly in the context
of front tracking methods for fluids and PDEs with smooth, fixed boundaries.
One of the earliest such methods is the immersed boundary method (IBM) of
Peskin, in which the boundary exerts an effect on a fluid, represented on a
rectangular mesh, using approximations to delta functions located on the
nonconforming smooth boundary (see~\cite{peskin1972flow}).
This method was extended by Leveque and Li to the immersed interface method
(IIM), which modifies the finite difference stencil in the vicinity of the
boundary to correct for error terms derived from the underlying Taylor
expansions (see~\cite{leveque1994immersed}).
In the explicit jump immersed interface method (EJIIM) proposed by Wiegmann
and Bube, the corrections of the IIM are applied directly in terms of the jumps
in the solution and its derivatives.  Importantly, when the jumps are known
a priori, the corrections are applied to the right-hand side of the
discretized system of
equations; when they are unknown, they are simultaneously solved for and
used to correct the solution in the same spirit as the IIM
(see~\cite{wiegmann2000explicit}). We note that the idea of applying
corrections to the right-hand side of the system was also suggested earlier
by Fornberg and Meyer-Spasche in~\cite{fornberg1992finite}, in which they
proposed a method for eliminating the first term in the expansion of the
error near the nonconforming boundary.  The ghost cell method (GSM) proposed
by Gibou, Fedkiw, and others \cite{gibou2002second}, based on the ghost fluid method (GFM) \cite{fedkiw1999non}, is an alternative way of applying the jump
corrections, in which ghost points are defined near the boundary, and
equations for their values are adjoined to the discretized system.
In \cite{gibou2002second}, the authors observe that a second-order
scheme can be constructed in which the discretized system is symmetric,
however, they also observe that the resulting finite difference matrices
becomes nonsymmetric for orders higher than two (see, for
example,~\cite{gibou2005fourth}).
%
%
All of the aforementioned methods assume that the jumps at the nonconforming
interface are either known beforehand, or are determined by augmenting the
finite difference system with additional equations.


We present a method that does not augment or alter the finite difference matrix, and does not assume that the jumps are known in advance.
We describe a high order deferred correction type algorithm for
computing both the solution and the free boundary of an LCP. The idea is to
derive the correction from the solution itself, after it has already been
computed without any correction, or with a correction of a lower order.  The
correction is then applied to the right-hand side, and the problem is
re-solved with the same matrix to one order of accuracy higher than before.  Two key ideas which
we use to rigorously justify this procedure are the smoothness of the error
away from the free boundary, which justifies the numerical differentiation
and extrapolation of the solution to obtain the jumps, and the fact that the
Green's function describing the error near the free boundary decreases like
$O(h)$ as the gridsize $h$ goes to zero, which is needed to show that the
jump corrections are computed to a sufficiently high order. Since the
corrections are computed separately and are applied exclusively to the right
hand side, the matrix of the system to be solved is identical to the
original finite difference matrix at each correction stage.
In fact, since the
solution at the previous correction stage can be used as an initial guess to
penalty iteration at the subsequent stage, only one or two iterations are
required for all correction stages after the first.  The jump corrections
are computed to high order by one-sided finite differences and
extrapolation, and the location of the free boundary is determined, also to
high order, from the solution by a combination of Lagrange interpolation and
Newton's method. The deferred correction procedure can, at least in
principle, be continued to indefinitely high orders, although we only apply it
to fourth-order. We also note that the principles behind our deferred
correction method are completely general, in the sense that they could be
applied to essentially any free boundary problem formulated as an LCP.  We
demonstrate the effectiveness of the method on several examples with a
one-dimensional space component, with and without a time component.  We also
apply our method to the well-known problem of pricing American put
options.

\section{Preliminaries}\label{sec:prelim}
\subsection{The LCP formulation of free and moving boundary problems}\label{sec:model-problem}
One form of the variational inequality representation of a free boundary problem in one dimension is
\begin{equation}\label{eqn:IVP-LCP}
\begin{cases}
&\partial_t\hat{V} - \mathcal L\hat{V} - g \geq 0,\\
&\hat{V} - V^* \geq 0,\\
&(\partial_t\hat{V} - \mathcal L\hat{V} - g)\cdot(\hat{V} - V^*) = 0,
\end{cases}
\end{equation}
see, for example \cite{friedman}, where $\mathcal L$ is a second-order differential operator
\begin{equation}\label{eqn:L-operator}
    \mathcal L = p(t,S)\frac{\partial^2}{\partial S^2} + w(t,S)\frac{\partial}{\partial S} + z(t,S),
\end{equation}
$V^*(t,S)$ is a given function, sometimes called the obstacle function or the payoff function, and $p$, $w$, $z$ and $g = g(t,S)$ are also given functions. Problem (\ref{eqn:IVP-LCP}) is also called a linear complementarity problem. Note that all three relations in \autoref{eqn:IVP-LCP} need to be satisfied. The solution of \autoref{eqn:IVP-LCP} is separated  into two parts by a moving boundary $S_f(t)$. The goal is to find the solution $\hat{V} = \hat{V}(t,S)$ such that either $\hat{V}-V^*>0$ and $\partial_t\hat{V} - \mathcal L\hat{V} - g = 0$, on what we call the PDE region of the solution, or $\partial_t\hat{V} - \mathcal L\hat{V} - g \geq 0$ and $\hat{V}-V^* = 0$, on what we call the penalty region of the solution.

In elliptic obstacle problems, the $\partial_t$ term disappears in the above formulation, and $\mathcal L$ is an elliptic operator. That is, the problem becomes
\begin{equation}\label{eqn:BVP-LCP}
\begin{cases}
&-\mathcal L\hat{V} - g \geq 0,\\
&\hat{V} - V^* \geq 0,\\
&( - \mathcal L\hat{V} - g)\cdot(\hat{V} - V^*) = 0,
\end{cases}
\end{equation}
In the American option pricing problems, $\partial_t -\mathcal L$ is the famous Black-Scholes operator with $\mathcal L = \mathcal L_{BS}$ and
\begin{equation}\label{eqn:BS-operator}
    \mathcal L_{BS} \equiv \frac{\sigma^2S^2}{2}\partial_{SS} + (r-d)S\partial_S - r,
\end{equation}
where $S$ is the underlying asset price, $r$ is the risk-free rate, $d$ is the dividend rate of the underlying asset, $\sigma$ is the volatility, and $t$ is the backward time from expiry. Typical payoff functions are
\begin{equation*}
    V^*(S) = \max\{S-K,0\}\; \text{ or } V^*(S) = \max\{K-S,0\}
\end{equation*}
for the American call and put options, respectively. Note that for American put and call options, the obstacle function is not time dependent.
It can be shown that the solution is only piecewise smooth, and the value matching and smooth pasting conditions
\begin{equation}\label{eqn:value-deriv-matching}
    \hat{V}(t, S_f(t)) = V^*(S_f(t)), \quad \frac{\partial \hat{V}}{\partial S}(t, S_f(t)) = \frac{\partial V^*}{\partial S}(S_f(t)),
\end{equation}
hold at the moving boundary (see, for example \cite{wilmott1995mathematics}), while the second derivative is discontiuous at $S_f(t)$. We see that the solution is only $C^1$ in space.

There has been much work that tries to achieve high-order convergence for solving the American option pricing problems.
Nonuniform-mesh techniques have been proposed to deal with the nonsmoothness in the solution at the undetermined exercise boundary \cite{oosterlee2005accurate,dilloo2017high}.
In the work by Oosterlee and Leentvaar \cite{oosterlee2005accurate}, the authors propose to use fourth-order finite differences in space and BDF4 in time, together with time-dependent grid-stretching in a predictor-corrector type scheme to achieve fourth-order accuracy. However, the authors of \cite{oosterlee2005accurate} did not provide numerical results on the convergence order for American options.
In this paper, we develop a general deferred correction algorithm using fourth-order finite difference method in space and BDF4 in time for solving free and moving boundary problems.

\subsection{Penalty method for solving the LCP}
In this paper, we solve the LCP using the penalty method as discussed in \cite{forsyth2002quadratic}. We approximate \autoref{eqn:IVP-LCP} by the penalized nonlinear PDE
\begin{equation}\label{eqn:IVP-penalty-LCP}
    \partial_t V = \mathcal LV + g + \rho \max\{V^* - V, 0\},
\end{equation}
for moving boundary problems, and
\begin{equation}\label{eqn:BVP-penalty-LCP}
    \mathcal LV + g + \rho \max\{V^* - V, 0\} = 0,
\end{equation}
for free boundary problems, where $\rho$ is a large positive penalty parameter, $\mathcal L$ is defined in \autoref{eqn:BS-operator}, and $V^*$ is the payoff function, as defined in Section \ref{sec:model-problem}, which also serves as the initial condition for PDE in (\ref{eqn:IVP-penalty-LCP}). When $\rho\rightarrow\infty$, either $V - V^* \geq 0$ or $V^* = V + \epsilon$ for $0 < \epsilon \ll 1$, where $\epsilon = \mathcal O(\rho^{-1})$ and $\rho \max\{V^* - V, 0\}$ is bounded, see \cite{friedman}. Using a finite volume discretization and applying the generalized Newton's iteration, also referred to as discrete penalty iteration, to the discretized PDE, the authors of \cite{forsyth2002quadratic} are able to prove monotonic convergence and finite termination of the algorithm under certain conditions. Moreover, second-order convergence can be obtained with an adaptive time step selector.

\section{Discretization, jump corrections and error analysis}\label{sec:appmaratus}
\subsection{Discretization of the penalized equation}\label{sec:discretization}
In this section, we describe the discretization of \autoref{eqn:BVP-penalty-LCP} and \autoref{eqn:IVP-penalty-LCP}, which will later lead to the formulation of a penalty iteration method for  solving \autoref{eqn:BVP-penalty-LCP} and \autoref{eqn:IVP-penalty-LCP}, similar to the second-order penalty method introduced in \cite{forsyth2002quadratic}. Unlike \cite{forsyth2002quadratic}, we use fourth-order finite difference space discretization and BDF4 time-stepping in order to obtain high-order accuracy. 

Consider a discretized domain $S_0 < S_1 < \cdots < S_{M+1}$ where $S_0$ and $S_{M+1}$ represent the left and right boundary respectively. Let $\tilde{V}_j^n \approx V(t_n, S_j)$ be the finite difference approximation to the true solution $V(t,S)$ of \autoref{eqn:IVP-penalty-LCP} at time $t_n$, and space point $S_j$. We drop the superscript $n$ when time is irrelevant. On a uniform grid with grid size $h$, the fourth-order finite difference approximation to $\frac{\partial^2V}{\partial S^2}(t, S_j)$ is given by the operator
\begin{equation*}
    D^2_4V_j \equiv \frac{1}{12h^2}(-V_{j-2} + 16V_{j-1} - 30V_j + 16V_{j+1} - V_{j+2}),
\end{equation*}
for $2 \leq j \leq M-1$, and
\begin{equation*}
    D^2_4V_1 \equiv \frac{1}{12h^2}(10V_0 - 15V_1 - 4V_2 + 14V_3 - 6V_4 + V_5),
\end{equation*}
for $j = 1$, and similarly for $j = M$.
On a nonuniform grid, the finite difference weights can be obtained by the method of undetertermined coefficients in a stable way (see, for example, \cite{fornberg1998classroom}). We denote the generic fourth-order finite difference approximation to $\frac{\partial^2V}{\partial S^2}$ at point $S_j$ to be
\begin{equation*}
    D_4^2V_j \equiv c_{-2}V_{j-2} + c_{-1}V_{j-1} + c_0V_j + c_{1}V_{j+1} + c_2V_{j+2},
\end{equation*}
where we abuse notation here and denote the finite difference coefficients at the points $x_{j-2},\; x_{j-1},\; x_j,$ $x_{j+1},\; x_{j+2}$ by $c_{-2},\; c_{-1},\; c_0,\; c_1,\;c_2$, respectively, for the finite difference approximation at $x_j$. Fourth-order finite difference discretization of the first derivative $\frac{\partial V}{\partial S}(t,S_j)$ can be obtained similarly using a five-point stencil, which we omit for brevity.

Let $\mathbf S$ denote the vector of the interior grid points, i.e. $\mathbf S = [S_1, \ldots, S_M]^T$.
Assuming Dirichlet boundary conditions, the fourth-order finite differences above give us the space discretization of $\partial_{SS}V$ and $\partial_SV$
\begin{equation*}
    \frac{\partial V}{\partial S}(t,\mathbf S) \approx \bar{\mathbf L}_1\tilde{\mathbf V}_{\text{aug}}, \quad \frac{\partial^2 V}{\partial S^2}(t,\mathbf S) \approx \bar{\mathbf L}_2\tilde{\mathbf V}_{\text{aug}},
\end{equation*}
where $\tilde{\mathbf V}_{\text{aug}} \equiv [\tilde{V}_0,\; \tilde{V}_1,\; \ldots,\; \tilde{V}_{M+1}]^T$ is the finite difference solution vector, $\bar{\mathbf L}_1$ and $\bar{\mathbf L}_2$ are $M\times(M+2)$ matrices with the coefficients of the corresponding finite difference stencil on each row. Let $\mathbf L$ be an $M\times M$ matrix defined by
\begin{equation}\label{eqn:L-mat-def}
    \mathbf L \equiv \mathbf P\mathbf L_2 + \mathbf W\mathbf L_1 + \mathbf Z,
\end{equation}
where $\mathbf L_2$ and $\mathbf L_1$ are $M\times M$ matrices from the interior columns of $\bar{\mathbf L}_2$ and $\bar{\mathbf L}_1$, respectively, and $\mathbf P$, $\mathbf W$, and $\mathbf Z$ are diagonal matrices with diagonal entries $[\mathbf P]_{jj}=p(t,S_j)$, $[\mathbf W]_{jj} = w(t,S_j)$, and $[\mathbf Z]_{jj} = z(t,S_j)$ for $j = 1, \ldots, M$.
Then the discretization of $\mathcal LV + g$ becomes
\begin{equation*}
    \mathcal LV(t,\mathbf S) + g(t,\mathbf S)\approx \mathbf L\tilde{\mathbf V} + \mathbf b,
\end{equation*}
where $\tilde{\mathbf V} \equiv [\tilde{V}_1,\; \tilde{V}_2,\; \ldots,\; \tilde{V}_{M}]^T$, and
\begin{align*}
    \mathbf b = &p(t,S_0)V(t,S_0)\bar{\mathbf L}_2[:,1] + w(t,S_0)V(t,S_0)\bar{\mathbf L}_1[:,1] \\
    &+ p(t,S_{M+1})V(t,S_{M+1})\bar{\mathbf L}_2[:,M+2] + w(t,S_{M+1})V(t,S_{M+1})\bar{\mathbf L}_1[:,M+2] + g(t,\mathbf S),
\end{align*}
which is a vector that incorporates the boundary conditions, where $\bar{\mathbf L}_1[:,j]$ and $\bar{\mathbf L}_2[:,j]$ denote the $j$-th columns of $\bar{\mathbf L}_1$ and $\bar{\mathbf L}_2$, respectively. The penalty term in \autoref{eqn:IVP-penalty-LCP} and \autoref{eqn:BVP-penalty-LCP} can be discretized by
\begin{equation}\label{eqn:discrete-penalty}
    \mathbf q(\tilde{\mathbf V}) \equiv \rho\bm{\mathcal I}_{\tilde{\mathbf V}}(\mathbf V^* - \tilde{\mathbf V})
\end{equation}
where $\mathbf V^* = [V^*_{1},\; V^*_{2},\; \ldots,\; V^*_{M}]^T$ is the vector of the payoff function values on the grid points $S_1$ to $S_M$, and $\bm{\mathcal I}_{\tilde{\mathbf V}}$ is a diagonal matrix whose diagonal entries are
\begin{equation}\label{eqn:penalty_identity}
    [\bm{\mathcal I}_{\tilde{\mathbf V}}]_{i,i} =
    \begin{cases}
    1, & V^*_{i} > \tilde{V}_i,\\
    0, & \text{else}.
    \end{cases}
\end{equation}
Therefore, we obtain the discretization of the right-hand side of \autoref{eqn:IVP-penalty-LCP},
\begin{equation}\label{eqn:discrete_penalty}
    \mathcal LV(t,\mathbf S) + g(t,\mathbf S) + \rho\max\{ V^*(\mathbf S)-V(t,\mathbf S), 0 \} \approx \mathbf L\tilde{\mathbf V} + \mathbf b + \mathbf q(\tilde{\mathbf V}).
\end{equation}
Assuming BDF4 uniform time discretization, and defining
\begin{equation}\label{eqn:A-mat-def}
   \mathbf A \equiv \frac{25}{12}\mathbf I - k\mathbf L,
\end{equation}
the complete discretization of \autoref{eqn:IVP-penalty-LCP} including time stepping follows the rule
\begin{equation}\label{eqn:BDF4-discrete}
    \mathbf A\tilde{\mathbf V}^{n+4} = 4\tilde{\mathbf V}^{n+3} - 3\tilde{\mathbf V}^{n+2} + \frac{4}{3}\tilde{\mathbf V}^{n+1} - \frac{1}{4}\tilde{\mathbf V}^{n} + k\mathbf b^{n+4} + k\mathbf q(\tilde{\mathbf V}^{n+4}),
\end{equation}
where $k$ is time step size, $\mathbf I$ is the identity matrix of size $M\times M$, and the superscript $n$ means the $n$-th time step. 
We also obtain the discretization of \autoref{eqn:BVP-penalty-LCP} as
\begin{equation}\label{eqn:BVP-discrete}
\mathbf L\tilde{\mathbf V} + \mathbf b + \mathbf q(\tilde{\mathbf V}) = 0.
\end{equation}
Systems (\ref{eqn:BDF4-discrete}) and \autoref{eqn:BVP-discrete} are nonlinear systems due to the presence of the penalty term, and we solve them using the penalty iteration described in \cite{forsyth2002quadratic}.


We note that, at this point, we do not specify the choice
of $\rho$ in the terms $\mathbf q(\tilde{\mathbf V}^{n+4})$ and $\mathbf q(\tilde{\mathbf V})$
in (\ref{eqn:BDF4-discrete}) and \autoref{eqn:BVP-discrete}, respectively.
This will be discussed in Subsection \ref{sec:alg-mb}.

\subsection{Finite difference approximation on a nonsmooth but piecewise smooth function}\label{sec:FD-nonsmoothness}
As has been mentioned in the previous section, the solution of the LCP has a discontinuous second derivative at the free boundary at all times. This is a major factor that causes the degeneracy of convergence rate when using the finite difference method on uniform grids. Since the finite difference approximation is based upon Taylor expansions, a certain level of smoothness has to be assumed in order to obtain the corresponding accuracy. When this smoothness requirement is not satisfied even at a single point, the truncation error will be contaminated by an additional error, and propagated to other points in the solution through the Green's function, as we will see later.
The analysis of this section is similar to the analysis of Li \cite{li1998fast}, and Wiegmann and Bube \cite{wiegmann2000explicit}, except that it is applied to our particular high-order finite difference operator.

\begin{figure}[h]
    \centering
\begin{tabular}{c c c}
    \includegraphics[width=0.32\textwidth]{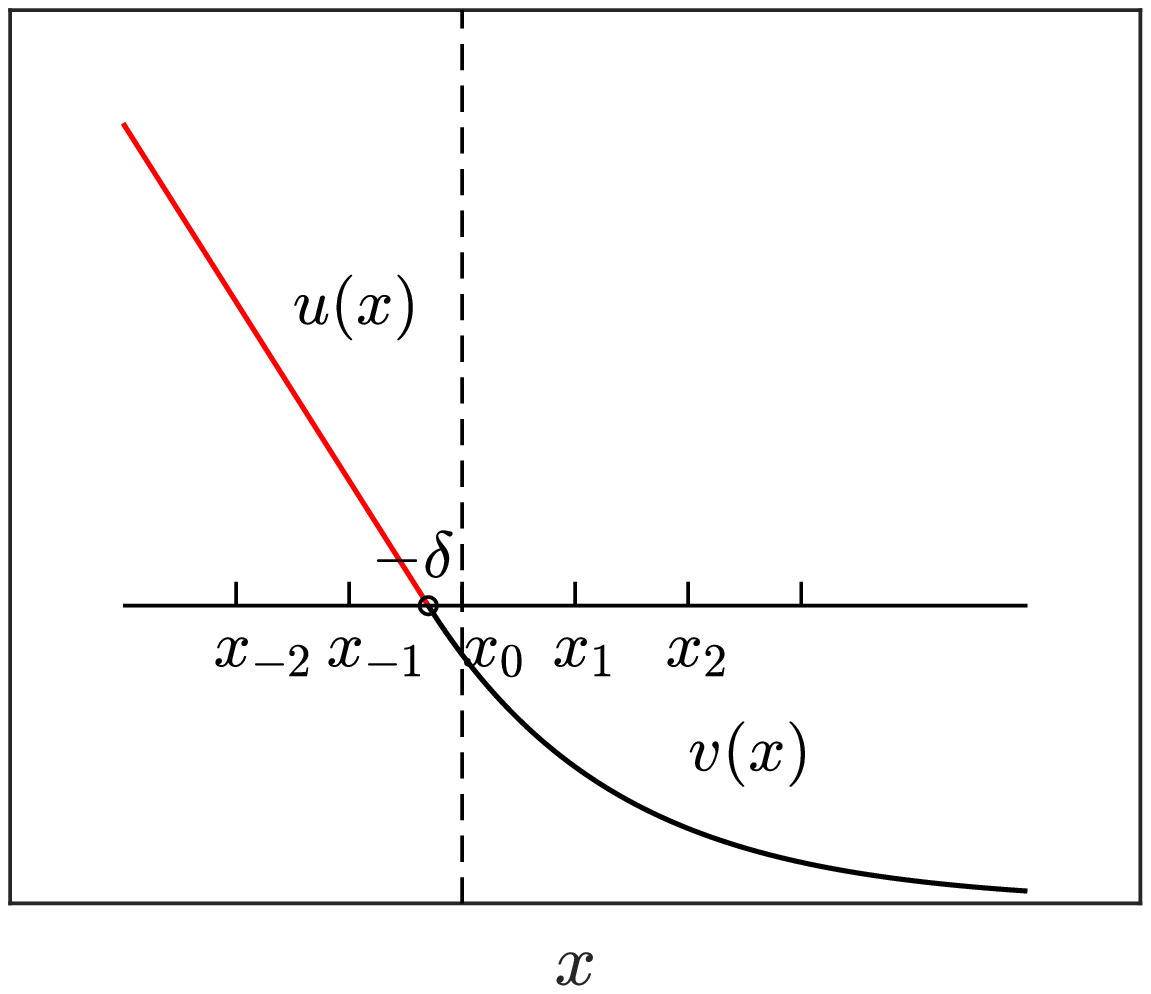} &  \includegraphics[width=0.32\textwidth]{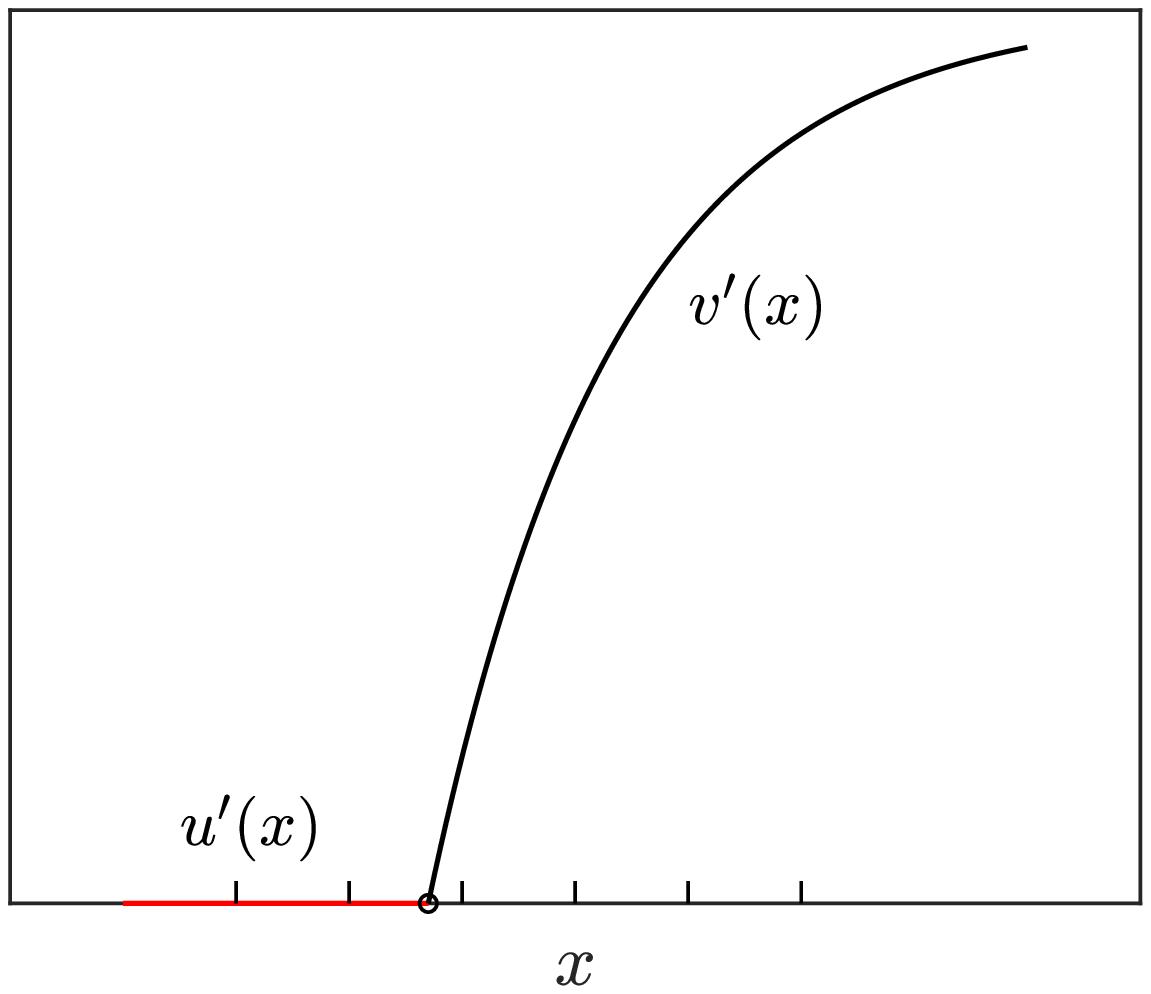} & \includegraphics[width=0.32\textwidth]{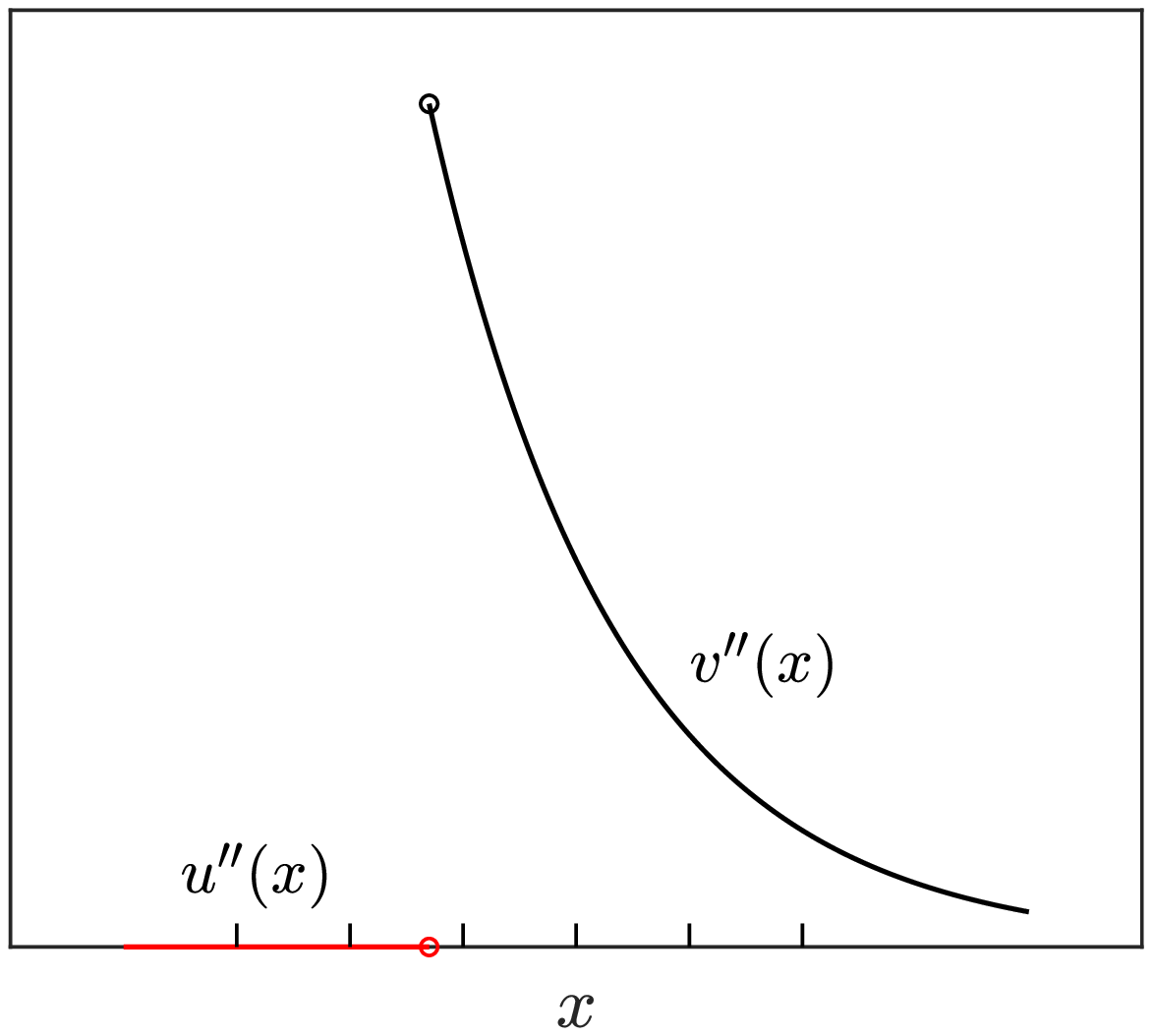}\\
    (a) $f(x)$ & (b) $f'(x)$ & (c) $f''(x)$
\end{tabular}
    \caption{An example graph of a nonsmooth function with a point of discontinuity of the second derivative at $x=-\delta$}
    \label{fig:nonsmooth-func-ex}
\end{figure}

To analyze the impact of piecewise smoothness on the finite difference approximation, consider a piecewise smooth function
\begin{equation}\label{eqn:nonsmooth-fx}
    f(x) =
    \begin{cases}
    v(x), & x+\delta > 0,\\
    u(x), & x+\delta \leq 0,
    \end{cases}
\end{equation}
such that $u(-\delta) = v(-\delta)$ and $u'(-\delta) = v'(-\delta)$, where $\delta$ is a positive constant. In addition, suppose that $u(x)$ and $v(x)$ admit smooth extensions, i.e., $u(x)$ is well defined and can be smoothly extended to the domain $x > -\delta$, and similarly $v(x)$ can be smoothly extended to $x < -\delta$. Let $\{x_j\}$ be a grid with $x_i < x_j$ for $i < j$, and with $x_{-1} < -\delta < x_0 = 0$. An example graph of function $f(x)$ with grid points $x_{-2}$ to $x_2$ is shown in Figure \ref{fig:nonsmooth-func-ex}. We want to approximate the second derivative of $f(x)$ at grid points around the nonsmooth position $x = -\delta$.

\subsubsection{Second-order finite difference scheme}
For notational convenience, we give a detailed derivation only for the second-order method. Derivations for the fourth-order method follow similarly. We pick the grid and location of the nonsmooth point only for the ease of demonstration. The following derivation is generalizable to any other function of the same form, irrespectively of where the nonsmooth point is located.
When using a second-order finite difference method to approximate the second derivative of $f(x)$ at point $x=x_0 = 0$, we are actually computing
\begin{equation}\label{eqn:second-center-fd-fake}
    D^2f_0 = \frac{1}{\bar{h}_{0}}\left[ \frac{1}{h_{0}}u_{-1} - \left( \frac{1}{h_{0}} + \frac{1}{h_{1}}\right)v_0 + \frac{1}{h_{1}}v_1 \right],
\end{equation}
where $u_j$, $v_j$ denote $u(x_j), v(x_j)$ respectively, $D^2$ represents the standard centered three-point finite-difference operator, $h_j = x_j - x_{j-1}$, and $\bar{h}_j = (h_j+h_{j+1})/2$. Note that the value of $u(x_{-1})$ instead of $v(x_{-1})$ is used for the left-most stencil point in \autoref{eqn:second-center-fd-fake}. This is because the finite difference operator is applied to $f(x)$, which is equal to $u(x_{-1})$ at point $x_{-1}$. However, the correct (in the sense that it is second-order accurate) approximation to the second derivative at the point $x_0$ should be
\begin{equation}\label{eqn:second-center-fd-real}
    D^2v_0 = \frac{1}{\bar{h}_{0}}\left[ \frac{1}{h_{0}}v_{-1} - \left( \frac{1}{h_{0}} + \frac{1}{h_{1}}\right)v_0 + \frac{1}{h_{1}}v_1 \right],
\end{equation}
where we recall the assumption that $v(s)$ has smooth extension for $x < -\delta$. Note that $u_{-1}$ in the formula $D^2f_0$ is replaced by $v_{-1}$ in the formula $D^2v_0$. The other problematic point is at $x = x_{-1} = -h_{0}$, where we approximate the derivative by
\begin{align*}
    D^2f_{-1} &= \frac{1}{\bar{h}_{-1}}\left[ \frac{1}{h_{-1}}u_{-2} - \left( \frac{1}{h_{-1}} + \frac{1}{h_{0}}\right)u_{-1} + \frac{1}{h_{0}}v_0 \right],
\end{align*}
rather than the second-order accurate finite difference
\begin{align*}
    D^2u_{-1} &= \frac{1}{\bar{h}_{-1}}\left[ \frac{1}{h_{-1}}u_{-2} - \left( \frac{1}{h_{-1}} + \frac{1}{h_{0}}\right)u_{-1} + \frac{1}{h_{0}}u_0 \right].
\end{align*}
The points $x_{-1}$ and $x_0$ are the only problematic points for a second-order method. The degeneracy of the finite difference approximation accuracy comes from the inconsistency between the formulas for $D^2f_0$ and $D^2v_0$, and between $D^2f_{-1}$ and $D^2v_{-1}$.

The following theorem describes the relationship between $D^2f_0, \; D^2f_{-1}$ and $D^2v_0, \; D^2u_{-1}$, respectively, in terms of the jumps of $u(x)$ and $v(x)$ at point $x = -\delta$, and quantifies the degeneration of accuracy.

\begin{theorem}\label{thm:second-order-correction}
Suppose $f(x)$ is given by \autoref{eqn:nonsmooth-fx}, where $f(x) = v(x)$ for $x>-\delta$ and $f(x) = u(x)$ for $x \leq -\delta$, with $u(-\delta) = v(-\delta)$ and $u'(-\delta) = v'(-\delta)$, where $u(x)$ and $v(x)$ admit smooth extensions. Consider the functions on a grid $\{x_j\}$ with $x_i < x_j$ for $i < j$, and with $x_{-1} < -\delta < x_0 = 0$. Then, $D^2f_0, \; D^2f_{-1}$ and $D^2v_0, \; D^2u_{-1}$ satisfy the relations
\begin{equation}\label{eqn:second-deriv-r}
\begin{aligned}
    D^2v_0 = D^2f_0 - \frac{(h_0-\delta)^2}{h_0(h_0+h_1)}(u''_{\delta}-v''_{\delta}) &+ \frac{(h_0-\delta)^3}{3h_0(h_0+h_1)}(u'''_{\delta}-v'''_{\delta}) \\
    &- \frac{(h_0-\delta)^4}{12h_0(h_0+h_1)}(u''''_{\delta}-v''''_{\delta}) + \mathcal O(h^3),\\
\end{aligned}
\end{equation}
and
\begin{equation}\label{eqn:second-deriv-l}
\begin{aligned}
    D^2u_{-1} = D^2f_{-1} + \frac{\delta^2}{h_0(h_{-1}+h_0)}&(u''_{\delta}-v''_{\delta}) + \frac{\delta^3}{3h_0(h_{-1}+h_0)}(u'''_{\delta}-v'''_{\delta}) \\
    &+ \frac{\delta^4}{12h_0(h_{-1}+h_0)}(u''''_{\delta}-v''''_{\delta}) + \mathcal O(h^3),
\end{aligned}
\end{equation}
where $h = \max\{h_0, h_1\}$, and the subscript $\delta$ denotes the quantities at the nonsmooth point $x = -\delta$, e.g. $u_{\delta}'' = u''(-\delta)$.
\end{theorem}

\begin{proof}
Subtracting \autoref{eqn:second-center-fd-fake} from (\ref{eqn:second-center-fd-real}), we get
\begin{equation}\label{eqn:D2f-D2v-compact}
    D^2f_0 = D^2v_0 + \frac{2}{h_0(h_0+h_1)}(u_{-1}-v_{-1}).
\end{equation}
Applying Taylor expansions for functions $u(x)$ and $v(x)$ around $x = -\delta$, we have
\begin{align*}
    u_{-1} = u_{\delta} - (h_0-\delta)u'_{\delta} + \frac{(h_0-\delta)^2}{2}u''_{\delta} - \frac{(h_0-\delta)^3}{6}u'''_{\delta} + \frac{(h_0-\delta)^4}{24}u''''_{\delta} + \mathcal O((h_0-\delta)^5),\\
    v_{-1} = v_{\delta} - (h_0-\delta)v'_{\delta} + \frac{(h_0-\delta)^2}{2}v''_{\delta} - \frac{(h_0-\delta)^3}{6}v'''_{\delta} + \frac{(h_0-\delta)^4}{24}v''''_{\delta} + \mathcal O((h_0-\delta)^5),
\end{align*}
which gives
\begin{equation}\label{eqn:u-v-at-delta}
\begin{aligned}
    u_{-1} - v_{-1} = \frac{(h_0-\delta)^2}{2}(u''_{\delta}-v''_{\delta}) &- \frac{(h_0-\delta)^3}{6}(u'''_{\delta}-v'''_{\delta}) \\
    &+ \frac{(h_0-\delta)^4}{24}(u''''_{\delta}-v''''_{\delta}) + \mathcal O((h_0-\delta)^5),
\end{aligned}
\end{equation}
using the assumptions that $u_{\delta} = v_{\delta}$ and $u'_{\delta} = v'_{\delta}$. Substituting \autoref{eqn:u-v-at-delta} into \autoref{eqn:D2f-D2v-compact}, we get \autoref{eqn:second-deriv-r}.
Following a similar derivation, we get \autoref{eqn:second-deriv-l}.
\end{proof}

\subsubsection{Fourth-order finite difference scheme}\label{sec:fourth-fd-nonsmooth}
In the previous section, we use the second-order approximation as a convenient way to demonstrate the essential relations that lead to our method. In this paper, we focus on high-order methods. Following exactly the same derivation procedure, we can arrive at similar formulas for fourth-order methods. The main difference between the second-order and fourth-order FDs is that in the fourth-order FDs there are four problematic points, namely $x_{-2}, x_{-1}, x_0, x_1$, instead of just two. Let the finite difference coefficients at the points $x_{j-2},\; x_{j-1},\; x_j,\; x_{j+1},\; x_{j+2}$ be denoted by $c_{-2},\; c_{-1},\; c_0,\; c_1,\;c_2$, respectively, for the finite difference approximation at $x_j$. We give the following theorem for fourth-order discretization.

\begin{theorem}\label{thm:fourth-order-correction}
Under the same assumptions as in \autoref{thm:second-order-correction}, we have that $D_4^2u_{-2}, \; D_4^2u_{-1}, \; D_4^2v_0, \; D_4^2v_{1}$ satisfy the relations
\begin{alignat}{2}
&D_4^2u_{-2} = D_4^2f_{-2} &&+ c_2\frac{\delta^2}{2}(u''_{\delta}-v''_{\delta}) + c_2\frac{\delta^3}{6}(u'''_{\delta}-v'''_{\delta}) + c_2\frac{\delta^4}{24}(u''''_{\delta}-v''''_{\delta}) + \mathcal O(h^3),\label{eqn:correction-1}\\
&D_4^2u_{-1} = D_4^2f_{-1} &&+ \left( c_1\frac{\delta^2}{2} + c_2\frac{(h_1+\delta)^2}{2} \right)(u''_{\delta}-v''_{\delta})\label{eqn:correction-2}\\
& &&+ \left( c_1\frac{\delta^3}{6}+c_2\frac{(h_1+\delta)^3}{6} \right)(u'''_{\delta}-v'''_{\delta}) \nonumber\\
& &&+ \left( c_1\frac{\delta^4}{24}+c_2\frac{(h_1+\delta)^4}{24} \right)(u''''_{\delta}-v''''_{\delta}) + \mathcal O(h^3),\nonumber\\
&D_4^2v_{0} = D_4^2f_{0} &&- \left( c_{-2}\frac{(h_{-1}+h_0-\delta)^2}{2} + c_{-1}\frac{(h_0-\delta)^2}{2} \right)(u''_{\delta}-v''_{\delta})\label{eqn:correction-3}\\
& &&+ \left( c_{-2}\frac{(h_{-1}+h_0-\delta)^3}{6}+c_{-1}\frac{(h_0-\delta)^3}{6} \right)(u'''_{\delta}-v'''_{\delta})\nonumber \\
& &&- \left( c_{-2}\frac{(h_{-1}+h_0-\delta)^4}{24} + c_{-1}\frac{(h_0-\delta)^4}{24} \right)(u''''_{\delta}-v''''_{\delta}) + \mathcal O(h^3),\nonumber\\
&D_4^2v_{1} = D_4^2f_{1} &&- c_{-2}\frac{(h_0-\delta)^2}{2}(u''_{\delta}-v''_{\delta}) + c_{-2}\frac{(h_0-\delta)^3}{6}(u'''_{\delta}-v'''_{\delta})\label{eqn:correction-4}\\
& &&- c_{-2}\frac{(h_0-\delta)^4}{24}(u''''_{\delta}-v''''_{\delta}) + \mathcal O(h^3),\nonumber
\end{alignat}
where $h = \max\{h_{-1}, h_0, h_1\}$.
\end{theorem}

\begin{proof}
From the approximation equations, we easily see that
\begin{align*}
&D_4^2f_{-2} = D_4^2u_{-2} + c_2(v_0-u_0),\\
&D_4^2f_{-1} = D_4^2u_{-1} + c_1(v_0-u_0) + c_2(v_1-u_1),\\
&D_4^2f_{0} = D_4^2v_{0} + c_{-2}(u_{-2}-v_{-2}) + c_{-1}(u_{-1}-v_{-1}),\\
&D_4^2f_{1} = D_4^2v_{1} + c_{-2}(u_{-1}-v_{-1}).
\end{align*}
Then, expressing the quantities $u_{-2}-v_{-2},\; u_{-1}-v_{-1},\;v_0-u_0,\;v_1-u_1$, by applying Taylor expansions to $u_{-2}, v_{-2}, u_{-1},v_{-1},v_0,u_0,v_1,u_1$ about the point $x = -\delta$, exactly as in the proof of \autoref{thm:second-order-correction}, we get the desired relations.
\end{proof}

From Theorems \ref{thm:second-order-correction} and \ref{thm:fourth-order-correction}, we see that, since $c_j = \mathcal O(1/h^2)$ and $\delta = \mathcal O(h)$, dominant $\mathcal O(1)$ terms appear in the truncation errors. Therefore, the second derivative approximations have degenerated orders of accuracy, regardless of the order of discretization.
In order to achieve the desired order of accuracy, we have to eliminate the remainder terms. This can be done by adding corrections. We call the right-hand side terms of Equations (\ref{eqn:correction-1})--(\ref{eqn:correction-4}) that are added to $D_4^2f_j$ the \emph{correction terms} to the finite difference approximation of the second derivatives at $x_{-2}$ to $x_1$. For example, we call $c_2\frac{\delta^2}{2}(u''_{\delta}-v''_{\delta}) + c_2\frac{\delta^3}{6}(u'''_{\delta}-v'''_{\delta}) + c_2\frac{\delta^4}{24}(u''''_{\delta}-v''''_{\delta})$ the correction terms to $D_4^2f_{-2}$ at $x_{-2}$. We will also refer to $u''_{\delta}-v''_{\delta}$, $u'''_{\delta}-v'''_{\delta}$ and $u''''_{\delta}-v''''_{\delta}$ as \emph{derivative jumps}.

\subsubsection{Modifying the finite differences with approximate corrections}\label{sec:fd-w-approx-correction}
The order of the FDs at the problematic points $x_{-2}$ to $x_1$ is determined by the dominant error term in the corrections. If we were able to apply the exact corrections using Equations \autoref{eqn:correction-1}--\autoref{eqn:correction-4}, we would fully recover the fourth-order convergence of the FDs at the four problematic points $x_{-2}$ to $x_1$. However, in this work, we assume that the exact free boundary location and the derivative jumps are not known a priori. They are instead approximated using a previously computed $\mathcal O(h^\ell)$ solution, as we describe later. Therefore, we replace the exact free boundary and derivative jumps in the correction terms by the approximated ones. The accuracies of the corrected FDs depend on the accuracy of the approximate free boundary and derivative jumps.

To see how the order of accuracy of the free boundary and derivative jumps affect the correction, suppose that the free boundary is known exactly. Then, it is obvious that $\mathcal O(h^\ell)$ derivative jumps will give rise to $\mathcal O(h^\ell)$ corrections.
On the other hand, suppose that the derivative jumps are known exactly, but we are given an approximate free boundary equal to $\delta+\mathcal O(h^\ell)$ with $1\leq \ell\leq 3$.
It is important to notice that the approximate free boundary introduces an extra source of error in the correction terms. To see this, we take one point, $x_{-2}$, for example.
The finite difference scheme with approximate correction terms becomes
\begin{equation}\label{eqn:sf-err-contrib}
\begin{aligned}
    &D_4^2f_{-2} + c_2\frac{(\delta+\mathcal O(h^{\ell}))^2}{2}(u_{\delta}''-v_{\delta}'') + c_2\frac{(\delta+\mathcal O(h^{\ell}))^3}{6}(u_{\delta}'''-v_{\delta}''') \\
    &\quad\quad\quad\quad\; + c_2\frac{(\delta+\mathcal O(h^{\ell}))^4}{24}(u_{\delta}''''-v_{\delta}'''')\\
    =& D_4^2f_{-2} + \left(c_2\frac{\delta^2}{2}+\mathcal O(h^{\ell-1})\right)(u_{\delta}''-v_{\delta}'') + \left(c_2\frac{\delta^3}{6}+\mathcal O(h^{\ell})\right)(u_{\delta}'''-v_{\delta}''') \\
    &\quad\quad\quad\quad +  \left(c_2\frac{\delta^4}{24}+\mathcal O(h^{\ell+1})\right)(u_{\delta}''''-v_{\delta}'''')\\
    =& D_4^2f_{-2} + c_2\frac{\delta^2}{2}(u_{\delta}''-v_{\delta}'') + c_2\frac{\delta^3}{6}(u_{\delta}'''-v_{\delta}''') + c_2\frac{\delta^4}{24}(u_{\delta}''''-v_{\delta}'''')\\
    &\quad\quad\quad\quad + \mathcal O(h^{\ell-1})(u_{\delta}''-v_{\delta}'') + \mathcal O(h^{\ell})(u_{\delta}'''-v_{\delta}''') + \mathcal O(h^{\ell+1})(u_{\delta}''''-v_{\delta}'''')\\
    =& D_4^2u_{-2} + \mathcal O(h^{\ell-1})(u_{\delta}''-v_{\delta}'') + \mathcal O(h^{\ell})(u_{\delta}'''-v_{\delta}''') + \mathcal O(h^{\ell+1})(u_{\delta}''''-v_{\delta}'''') + \mathcal O(h^3).
\end{aligned}
\end{equation}
Equation (\ref{eqn:sf-err-contrib}) implies that, when applying corrections using an approximate free boundary, the correction terms produce additional errors that are one order lower than the accuracy of the approximate free boundary. In order to improve the order of accuracy of the finite difference scheme by adding back the correction terms, we see that $\ell$ has to satisfy $\ell \geq 2$, because if $\ell = 1$, the leading order term of the corrections on the right-and side of \autoref{eqn:sf-err-contrib} is still of constant order $\mathcal O(h^{\ell-1})=\mathcal O(1)$. Therefore, we require the approximate derivative jumps and the free boundary location to be of at least $\mathcal O(h)$ and $\mathcal O(h^2)$, respectively, in order to increase the order of accuracy of the corrected finite differences to first-order, $\mathcal O(h^2)$ and $\mathcal O(h^3)$ to increase the order of accuracy to second-order, and so on.



\subsection{Convergence  of the fourth-order finite difference space discretization and its error propagation through the Green's function}\label{sec:error-and-greens-analysis}
\subsubsection{Boundary value problems}\label{subsec:bvp-conv}
In the previous section, we derived the correction terms for the finite difference approximations of derivatives of a nonsmooth but piecewise smooth function. Unless the values of $u''_{\delta}-v''_{\delta},\; u'''_{\delta}-v'''_{\delta},\;u''''_{\delta}-v''''_{\delta}$ are known, we cannot make use of these formulas directly to obtain fourth-order convergence. To solve this problem, we use a deferred correction approach, and successively compute the approximate derivatives from the lower-order solutions that are already known, and make sure to match up the orders of solutions and orders of corrections. In order to decide how much accuracy is required for the derivative approximation to result in corrections of the required order, we need to understand the error behaviour.

We consider boundary value problems that are time-independent, i.e., we consider \autoref{eqn:BVP-LCP}
so that we can leave the complexity of time evolution for later discussion. 
To analyze the error behavior of the space discretization scheme, we consider the finite difference approximation of the PDE in \autoref{eqn:BVP-penalty-LCP} given by \autoref{eqn:BVP-discrete}.



The theorem below describes the error behaviour of the fourth-order finite difference scheme applied to \autoref{eqn:BVP-penalty-LCP}, and how the nonsmoothness at the free boundary causes the convergence order of the fourth-order difference scheme to degenerate.

\begin{proposition}\label{prop:fourth-order-fd-error}
Consider the penalized PDE \autoref{eqn:BVP-penalty-LCP} with $V(S)$ being its exact solution, and the original LCP \autoref{eqn:BVP-LCP} with $\hat{V}(S)$ being its exact solution. Suppose that the first $m+1$ points $\hat{V}(S_0)$, $\hat{V}(S_1)$, \ldots, $\hat{V}(S_m)$ lie on the penalty region, i.e. $\hat{V}(S_j) = V^*(S_j) = V(S_j) \pm \epsilon$ for $0 \leq j \leq m$, with $0 < \epsilon \ll 1$ being approximately the size of the stopping tolerance set in the penalty iteration, and $\hat{V}(S_j) > V^*(S_j)$, $V(S_j)\pm\epsilon > V^*(S_j)$, for $m+1 \leq j \leq M+1$. Assume also that the approximate solution $\tilde{\mathbf V}$ of the penalty iteration exactly recovers $\bm{\mathcal I}_{\hat{\mathbf V}}$, i.e., $\bm{\mathcal I}_{\tilde{\mathbf V}} = \bm{\mathcal I}_{\hat{\mathbf V}}$.
Then, the error, $\mathbf e = [\hat{V}(S_1) - \tilde{V}_1,\; \hat{V}(S_2) - \tilde{V}_2,\; \ldots,\; \hat{V}(S_M) - \tilde{V}_M]^T$, of the fourth-order finite difference scheme in \autoref{eqn:BVP-discrete} for solving the penalized PDE \autoref{eqn:BVP-penalty-LCP} satisfies
\begin{equation}\label{eqn:error-decompose}
    (\mathbf L - \rho\bm{\mathcal I}_{\hat{\mathbf V}})\mathbf e =\bm{\upgamma} +
    \sum_{j=m-1}^{m+2}\mathcal O(1)\mathbf 1_j  + \sum_{j=1}^M \mathcal O(h^4)\mathbf 1_j
    \equiv \mathbf r,
\end{equation}
when the grid point $S_m$ is not exactly on the free boundary, i.e. $S_m < S_f$, where $\mathbf 1_j$ is the $j$-th column of an $M\times M$ identity matrix, and $[\bm{\upgamma}]_j =  (\mathcal L\hat{V}(S_j) + g(S_j))\mathbbm 1_{1\leq j \leq m}$, for $j = 1, \ldots, M$, where $\mathbbm 1_{1\leq j \leq m}$ is the indicator function,
which is one when $1\leq j \leq m$ and zero otherwise.
When $S_m$ is exactly on the free boundary, i.e., $S_m = S_f$, the sum in the second summation term is taken from $j = m-1$ to $m+1$.
\end{proposition}

\begin{proof}
Since $S_0, S_1, \ldots, S_m$ lie on the penalty region and we assume Dirichlet boundary conditions, $\bm{\mathcal I}_{\hat{\mathbf V}}$ is an $M\times M$ diagonal matrix with the diagonal elements $\left(\bm{\mathcal I}_{\hat{\mathbf V}}\right)_{i,i} = 1$ for $i = 1, \ldots, m$ and $\left(\bm{\mathcal I}_{\hat{\mathbf V}}\right)_{i,i} = 0$ for $i = m+1, \ldots, M$.
Hence, from the assumption that $\bm{\mathcal I}_{\tilde{\mathbf V}} = \bm{\mathcal I}_{\hat{\mathbf V}}$, we have
\begin{equation}\label{eqn:penalty2error}
\begin{aligned}
    \mathbf q(\tilde{\mathbf V}) &= \rho[V^*(S_1) - \tilde{V}_1,\; \ldots,\; V^*(S_m) - \tilde{V}_m,\; 0,\; \ldots, 0]^T\\
    &=\rho[\hat{V}(S_1) - \tilde{V}_1,\; \ldots,\; \hat{V}(S_m) - \tilde{V}_m,\; 0,\; \ldots, 0]^T\\
    &= \rho\bm{\mathcal I}_{\tilde{\mathbf V}}\mathbf e,
\end{aligned}
\end{equation}
Assume $S_m < S_f$, i.e. the grid point $S_m$ is not exactly on the free boundary. The proof for the case when $S_m = S_f$ is similar. 

From \autoref{thm:fourth-order-correction} for fourth-order discretization, we apply the discrete $\mathbf L$ operator to the true solution $\hat{\mathbf V} \equiv \hat{V}(\mathbf S)$ to get
\begin{equation}\label{eqn:fourth-order-BS-discretization2}
\begin{aligned}
    &\mathbf L\hat{\mathbf V} + \mathbf b
    = \sum_{j=1}^M\left(\mathcal L\hat{V}(S_j) + g(S_j) + \mathcal O(h^4)\right)\mathbf 1_j + \sum_{j=m-1}^{m+2}\mathcal O(1)\mathbf 1_j
    \equiv \mathcal L\hat{V}(\mathbf S) + g(\mathbf S) + \bm{\uptheta}, \\
    &\bm{\uptheta} \equiv  \sum_{j=m-1}^{m+2}\mathcal O(1)\mathbf 1_j + \sum_{j=1}^M \mathcal O(h^4)\mathbf 1_j.
\end{aligned}
\end{equation}
Since $\mathcal L\hat{V} + g = 0$ for $S_{m+1} \leq S \leq S_M$, and using \autoref{eqn:fourth-order-BS-discretization2},
we have
\begin{equation}\label{eqn:dt-L-exact-sub}
     \mathbf L\hat{\mathbf V} + \mathbf b = \mathcal L\hat{V}(\mathbf S) + g(\mathbf S) + \bm{\uptheta} = \bm{\upgamma} + \bm{\uptheta}.
\end{equation}
Subtracting \autoref{eqn:BVP-discrete} from \autoref{eqn:dt-L-exact-sub} and applying \autoref{eqn:penalty2error}, we get
\begin{equation*}
    \mathbf L(\hat{\mathbf V} - \tilde{\mathbf V}) - \mathbf q(\tilde{\mathbf V})  - \bm{\upgamma}  - \bm{\uptheta} = (\mathbf L - \rho\bm{\mathcal I}_{\hat{\mathbf V}})\mathbf e  - \bm{\upgamma} - \mathbf \bm{\uptheta} = \mathbf 0.
\end{equation*}
Therefore, the error satisfies
\begin{equation*}
    (\mathbf L - \rho\bm{\mathcal I}_{\hat{\mathbf V}})\mathbf e =
   \bm{\upgamma}+\sum_{j=m-1}^{m+2}\mathcal O(1)\mathbf 1_j +  \sum_{j=1}^M \mathbf O(h^4)\mathbf 1_j.
\end{equation*}
\end{proof}

\autoref{prop:fourth-order-fd-error} identifies the error equation $\mathbf e = (\mathbf L - \rho\bm{\mathcal I}_{\tilde{\mathbf V}})^{-1}\mathbf r$. The following proposition tells us how the operator $(\mathbf L - \rho\bm{\mathcal I}_{\tilde{\mathbf V}})^{-1}$ behaves.

\begin{proposition}\label{prop:penalized-matrix-inverse}
Consider the partitioning of the matrix $\mathbf L$
representing the discretization of \autoref{eqn:L-operator} and defined in \autoref{eqn:L-mat-def} as
\begin{equation}\label{eqn:L-mat}
    \mathbf L =
    \begin{bmatrix}
    \mathbf L_{11} & \mathbf L_{12}\\
    \mathbf L_{21} & \mathbf L_{22}
    \end{bmatrix},
\end{equation}
where
the submatrices
$\mathbf L_{11},\; \mathbf L_{12},\; \mathbf L_{21},\; \mathbf L_{22}$
are of sizes $m\times m$, $m\times(M-m)$, $(M-m)\times m$ and
$(M-m)\times(M-m)$ respectively, and $m$ is such that $S_m \le S_f$.
Assume $\mathbf L_{11}$ and $\mathbf L_{22}$ are nonsingular, and
$\rho$ is a positive number such that
$\rho \gg \max_{ij}\{|L_{i,j}|\}$. Assume also that
      $\max_{j}  \{|[\mathbf L_{22}^{-1}]_{1,j}|\} = \mathcal O(h^2)$,
      $\max_{j}  \{|[\mathbf L_{22}^{-1}]_{2,j}|\} = \mathcal O(h^2)$,
      $\max_{i}  \{|[\mathbf L_{22}^{-1}]_{i,1}|\} = \mathcal O(h^2)$,
and   $\max_{i}  \{|[\mathbf L_{22}^{-1}]_{i,2}|\} = \mathcal O(h^2)$.
Let $\bm{\mathcal I}$ be a diagonal matrix such that
$(\bm{\mathcal I})_{i,i} = 1$ for $i = 1,\ldots,m$,
$(\bm{\mathcal I})_{i,i} = 0$ for $i = m+1,\ldots,M$.
Then $(\mathbf L - \rho\bm{\mathcal I})^{-1}$ has the approximation
\begin{equation}
    (\mathbf L - \rho\bm{\mathcal I})^{-1}
     \approx
    \begin{bmatrix}
    \mathbf 0 & \mathbf 0\\
    \mathbf 0 & \mathbf L_{22}^{-1}
    \end{bmatrix}.
\end{equation}
\end{proposition}

\begin{remark1}
We denote $(\mathbf L_{22})^{-1}$ by $\mathbf L_{22}^{-1}$ for notation simplicity.
\end{remark1}

\begin{proof}
Since $\mathbf L_{11},\; \mathbf L_{22}$ are nonsingular,
the exact inverse matrix of $\mathbf L - \rho\bm{\mathcal I}$ is
\begin{equation*}
    (\mathbf L - \rho\bm{\mathcal I})^{-1} =
    \begin{bmatrix}
    \mathbf B & -\mathbf B \mathbf L_{12} \mathbf L_{22}^{-1}\\
    -\mathbf L_{22}^{-1} \mathbf L_{21} \mathbf B &
     \mathbf L_{22}^{-1} + \mathbf L_{22}^{-1} \mathbf L_{21} \mathbf B \mathbf
L_{12} \mathbf L_{22}^{-1}
    \end{bmatrix},
\end{equation*}
where
$\mathbf B = (\mathbf L_{11} - \rho\mathbf I
           - \mathbf L_{12}\mathbf L_{22}^{-1}\mathbf L_{21})^{-1}$,
and $\mathbf I$ is the identity matrix of size $m\times m$.
Note that $\mathbf L_{12}$ and $\mathbf L_{21}$
have only three nonzero entries in the lower-left and upper-right corners,
respectively, and all these entries are $\mathcal O(1/h^2)$.
The assumptions of the proposition
together with the special form of
$\mathbf L_{12}$ and $\mathbf L_{21}$
lead to $\mathbf L_{12}\mathbf L_{22}^{-1}\mathbf L_{21}$
have only four nonzero entries in the lower-right $2 \times 2$ corner
and these entries be $\mathcal O(1/h^2)$.
Since $\rho \gg \max_{ij}\{|L_{i,j}|\} = \mathcal O(1/h^2)$, we have
$\mathbf B \approx -\frac{1}{\rho}\mathbf I$. Therefore,
\begin{equation*}
    (\mathbf L - \rho\bm{\mathcal I})^{-1} \approx
    \begin{bmatrix}
    -\frac{1}{\rho}\mathbf I & \frac{1}{\rho}\mathbf L_{12} \mathbf L_{22}^{-1}\\
    \frac{1}{\rho}\mathbf L_{22}^{-1} \mathbf L_{21} &
    \mathbf L_{22}^{-1} - \frac{1}{\rho}
       \mathbf L_{22}^{-1} \mathbf L_{21} \mathbf L_{12} \mathbf L_{22}^{-1}
    \end{bmatrix} \approx
    \begin{bmatrix}
    \mathbf 0 & \mathbf 0\\
    \mathbf 0 & \mathbf L_{22}^{-1}
    \end{bmatrix},
\end{equation*}
where we have also taken into account that,
due to the special form of $\mathbf L_{21}$ and $\mathbf L_{12}$,
each entry of
$\mathbf L_{22}^{-1} \mathbf L_{21} \mathbf L_{12} \mathbf L_{22}^{-1}$
is composed as the sum of two terms of $\mathcal O(1)$.
\end{proof}

\begin{remark1}
Note that $\mathbf L_{22}$ behaves as the fourth-order
finite difference discretization of $\mathcal L$ on the grid
$S_{m+1},\; S_{m+2},\; \ldots,\; S_M$ with $S_m$ and $S_{M+1}$
as the boundary points. As such, the assumptions in the proposition that 
      $\max_{j}  \{|[\mathbf L_{22}^{-1}]_{1,j}|\} = \mathcal O(h^2)$,
      $\max_{j}  \{|[\mathbf L_{22}^{-1}]_{2,j}|\} = \mathcal O(h^2)$,
      $\max_{i}  \{|[\mathbf L_{22}^{-1}]_{i,1}|\} = \mathcal O(h^2)$,
and   $\max_{i}  \{|[\mathbf L_{22}^{-1}]_{i,2}|\} = \mathcal O(h^2)$ are typically true, under some conditions for the coefficients $p, w, z$ and the grid spacing. See, for example, \autoref{prop:greens} below.
\end{remark1}

\noindent
Let $\mathbf r$ be defined as in \autoref{prop:fourth-order-fd-error}
and $\mathbf r_{m+1:M}$ denote the subvector of $\mathbf r$ starting from
entry $m+1$ to $M$. We then have the following theorem.

\begin{theorem}\label{thm:penalty-error-equation}
Under the same assumptions as in \autoref{prop:penalized-matrix-inverse} and the assumption that $\bm{\mathcal I}_{\tilde{\mathbf V}} = \bm{\mathcal I}_{\hat{\mathbf V}}$,
when using the fourth-order finite difference scheme in
\autoref{eqn:BVP-discrete} to solve the penalized PDE
\autoref{eqn:BVP-penalty-LCP}, the error satisfies
\begin{equation*}
    \mathbf e \approx
    \begin{bmatrix}
    \mathbf 0\\
    \mathbf L_{22}^{-1}\mathbf r_{m+1:M}
    \end{bmatrix}.
\end{equation*}
\end{theorem}
\begin{proof}
The theorem is easily obtained from Propositions \ref{prop:fourth-order-fd-error}
and \ref{prop:penalized-matrix-inverse}.
\end{proof}

\autoref{thm:penalty-error-equation} shows that the penalty method
obtains the exact solution within a pre-specified tolerance
on the penalty region, while, on the PDE region,
where the solution satisfies the PDE \autoref{eqn:BVP-penalty-LCP},
the error is given by $\mathbf L_{22}^{-1}\mathbf r_{m+1:M}$.
Note also that,
on a uniform grid with step size $h$, $\mathbf L_{22}^{-1}$
can be thought of as the finite difference analogue
of the continuous Green's function of $\mathcal L$ on the PDE region,
scaled by $\frac{h}{2}$.
Therefore, we have
\begin{equation}\label{eqn:error-greens-prop}
    \mathbf e \approx \sum_{j=m+1}^M r_j \begin{bmatrix}
    \mathbf 0\\
    [\mathbf L_{22}^{-1}]_{:,j}
    \end{bmatrix} \approx \sum_{j=m+1}^M r_j \frac{h}{2}
    \begin{bmatrix}
    \mathbf 0\\
    G(\mathbf S_2, S_j)
    \end{bmatrix},
\end{equation}
where $G(\mathbf S_2, S_j)$ is a column vector of Green's function values $G(S,S_j)$ at points $\mathbf S_2 \equiv \{S_{m+1}, \ldots, S_M\}$. 

\begin{figure}[h]
\begin{tabular}{c c}
    \includegraphics[width=0.45\textwidth]{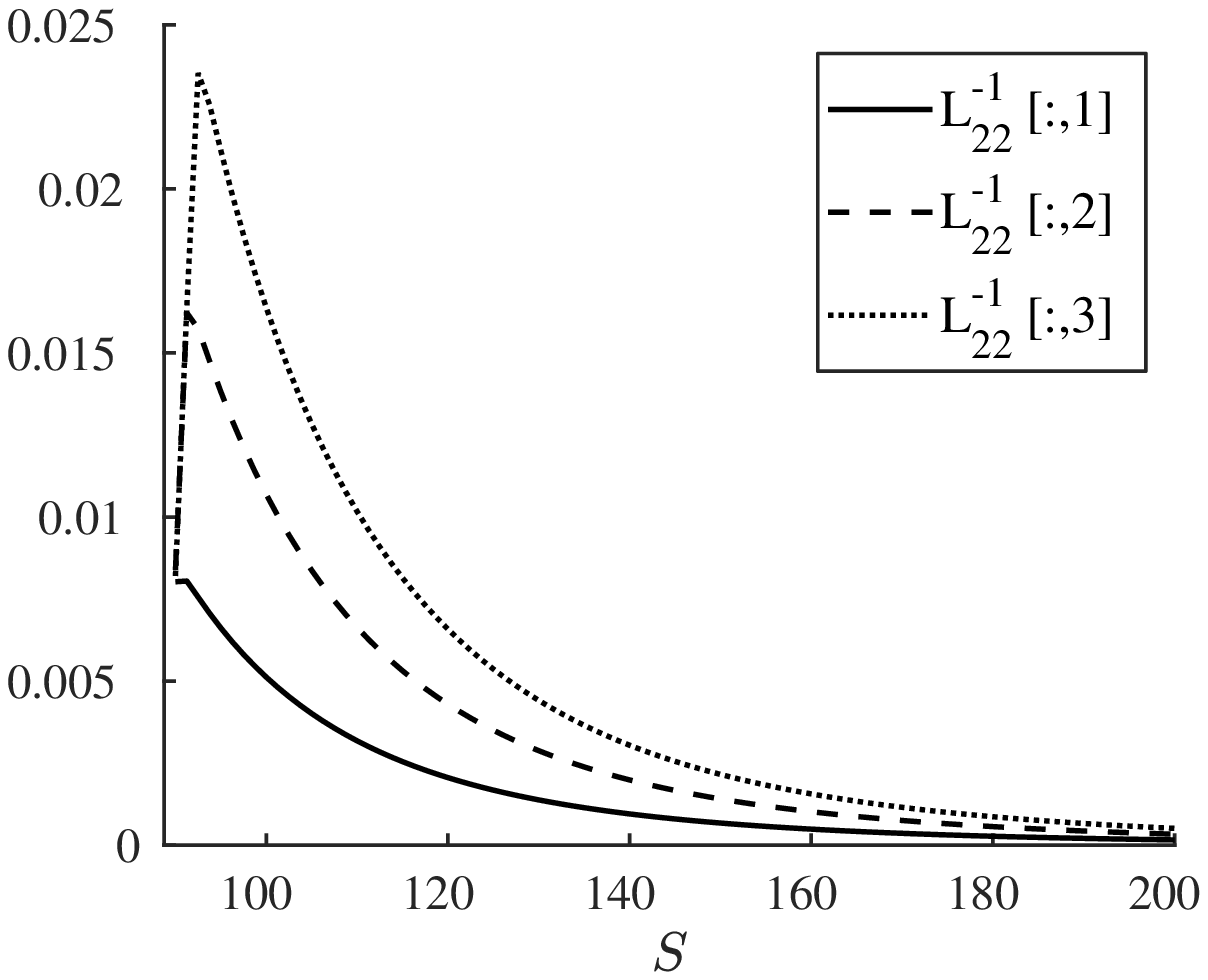} &  \includegraphics[width=0.45\textwidth]{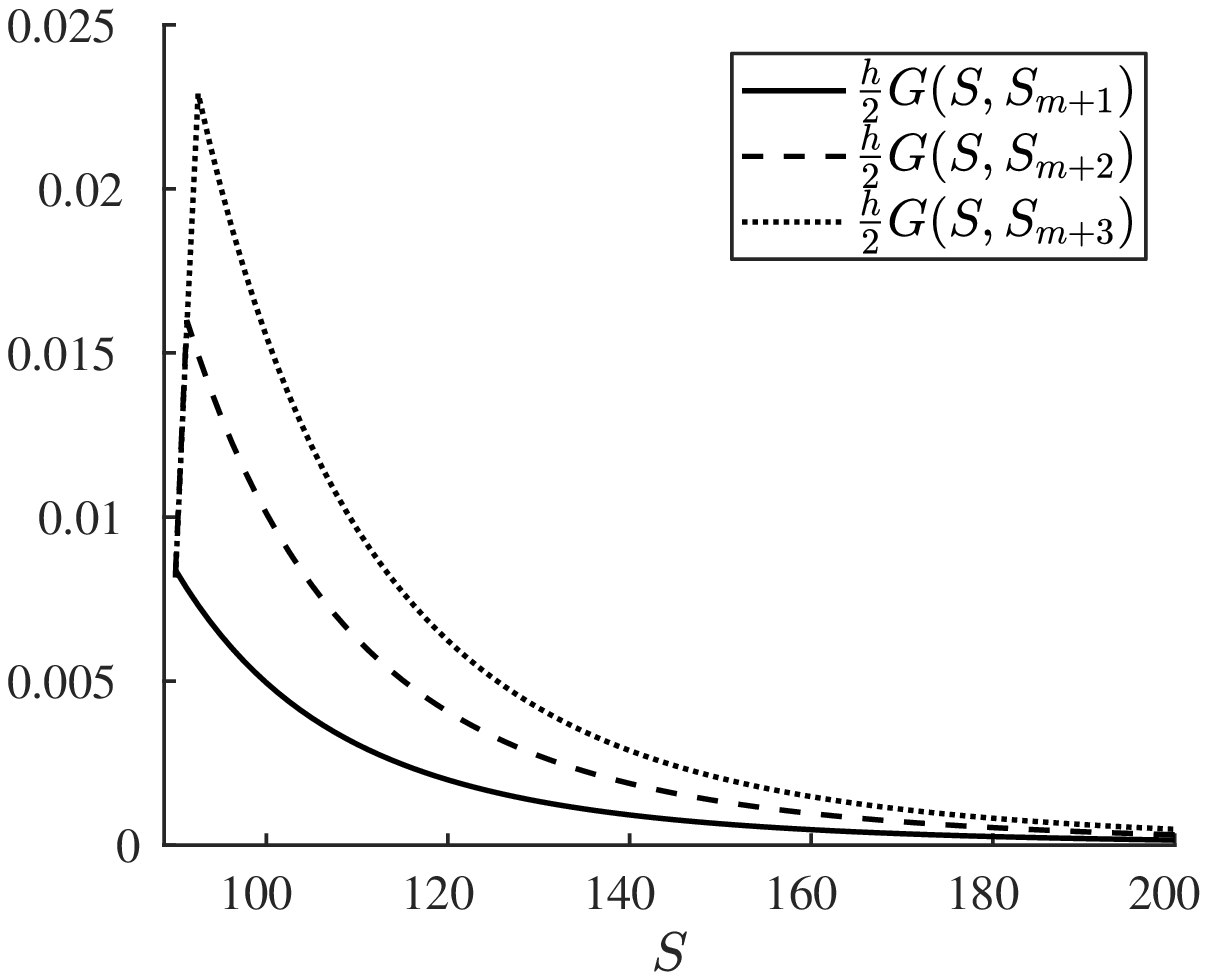}\\
    (a) The first 3 columns of $\mathbf L_{22}^{-1}$ & (b)
    The continuous Green's function $\frac{h}{2}G(S,S_j)$
\end{tabular}
    \centering
    \caption{(a) The first three columns of $\mathbf L_{22}^{-1}$ on an example uniform grid of size $h = 1.25$; (b) The scaled continuous Green's function $\frac{h}{2}G(S,S_j)$ for the operator $\mathcal L_{BS}$ at $S_{m+1}$, $S_{m+2}$, $S_{m+3}$ and for $S\geq S_{m+1}$. The free boundary location is $S_f = 89.748$.}
    \label{fig:L22inv}
\end{figure}

To visualize the size of entries of $\mathbf L_{22}^{-1}$ and its relation to the continuous Green's function, Figure \ref{fig:L22inv} gives the first three columns of the $\mathbf L_{22}^{-1}$ and corresponding $\frac{h}{2}G(S,S_{m+j})$, for $j = 1,\; 2$ and $3$, for the operator $\mathcal L_{BS}$ given by Equation (\ref{eqn:BS-operator}), on an example nonuniform grid where the free boundary is located at $S_m < S_f = 89.748 < S_{m+1}$. We can see that $[\mathbf L_{22}^{-1}]_{:,j}$ and $\frac{h}{2}G(S,S_{m+j})$ behave similarly.

In order to analyze the error behavior, we turn to understanding the properties of the Green's function $G(S,S_j)$, which is easier to investigate than its discrete analogue $\mathbf L_{22}^{-1}$. The following proposition gives the exact expression of the Green's function to a general operator.

\begin{proposition}\label{prop:greens}
Suppose that $Tu(x) = 0$ is a constant-coefficient, second-order homogeneous differential equation defined on the domain $[a, b]$. Let $u(x) = c_1e^{\xi_1 x} + c_2e^{\xi_2 x}$ denote the general solution to this equation. Suppose further that $\xi_1$ and $\xi_2$ are real and $\xi_1\neq\xi_2$. Then, the Green's function for the operator $T$ is
\begin{equation}\label{eqn:general-greens-function}
    G(x,\bar{x}) =
    \begin{cases}
    \dfrac{e^{(\xi_2-\xi_1)b}-e^{(\xi_2-\xi_1)\bar{x}}}{(\xi_2-\xi_1)e^{\xi_2\bar{x}}\left(e^{(\xi_2-\xi_1)b}-e^{(\xi_2-\xi_1)a}\right)}\left( e^{\xi_2x} - e^{(\xi_2-\xi_1)a+\xi_1x} \right), & a \leq x < \bar{x},\\
    \dfrac{e^{(\xi_2-\xi_1)a}-e^{(\xi_2-\xi_1)\bar{x}}}{(\xi_2-\xi_1)e^{\xi_2\bar{x}}\left(e^{(\xi_2-\xi_1)b}-e^{(\xi_2-\xi_1)a}\right)}\left( e^{\xi_2x} - e^{(\xi_2-\xi_1)b+\xi_1x} \right), & \bar{x} \leq x \leq b.
    \end{cases}
\end{equation}
Moreover, for any $x \in [a, b]$, we have
\begin{equation*}
    G(x,\bar{x}) = \mathcal O(\bar{x}-a), \quad \text{as } \bar{x}\rightarrow a.
\end{equation*}
\end{proposition}
\begin{proof}
The computation of the Green's function follows the standard procedure and we omit it. When $\bar{x}\leq x\leq b$, we have
\begin{equation*}
    e^{(\xi_2-\xi_1)a}-e^{(\xi_2-\xi_1)\bar{x}} = e^{(\xi_2-\xi_1)a}\left( 1 - e^{(\xi_2-\xi_1)(\bar{x}-a)} \right) \approx (\xi_1-\xi_2)(\bar{x}-a)e^{(\xi_2-\xi_1)a} = \mathcal O(\bar{x}-a),
\end{equation*}
as $\bar{x}\rightarrow a$. When $a\leq x\leq\bar{x}$, we have
\begin{equation*}
    e^{\xi_2x} - e^{(\xi_2-\xi_1)a+\xi_1x} = e^{\xi_2a}\left( e^{\xi_2(x-a)}\!-\!e^{\xi_1(x-a)} \right)\approx e^{\xi_2a}(\xi_2-\xi_1)(x-a)\leq e^{\xi_2a}(\xi_2-\xi_1)(\bar{x}-a)= \mathcal O(\bar{x}-a),
\end{equation*}
as $\bar{x}\rightarrow a$. Therefore, we see that $G(x,\bar{x}) = \mathcal O(\bar{x}-a)$ as $\bar{x}\rightarrow a$ for all $a\leq x\leq b$.
\end{proof}

\begin{figure}[h]
    \centering
    \includegraphics[width=0.5\textwidth]{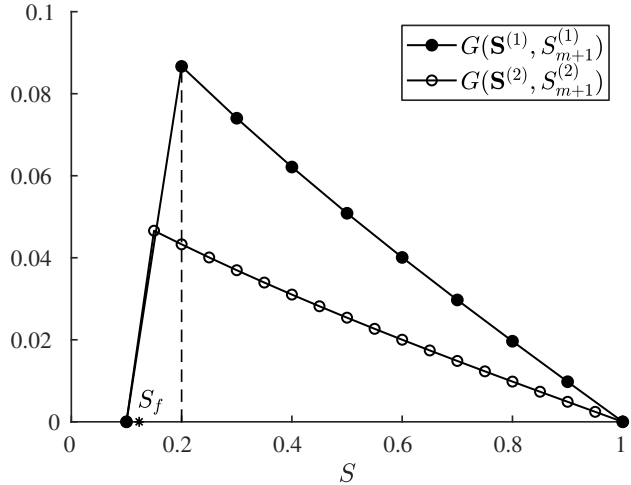}
    \caption{Illustration of Green's functions on two successive grid refinements $\mathbf S^{(1)}, \mathbf S^{(2)}$, for a hypothetical second-order differential equation.}
    \label{fig:bvp_greens}
\end{figure}

From Equation (\ref{eqn:error-greens-prop}), it is obvious that the error is dominated by the $\mathcal O(1)$ entries $r_{m+1}$ and $r_{m+2}$ (if $S_m$ is not on the free boundary) with the respective terms being propagated to every other grid point through the discrete Green's function.
Hence, at first glance, it appears that solving \autoref{eqn:BVP-discrete} using fourth-order finite differences will only give us first-order $\mathcal O(h)$ convergence.
However, from Proposition \ref{prop:greens}, with $a = S_f$, and from Equation \autoref{eqn:error-greens-prop}, we observe that second-order $\mathcal O(h^2)$ convergence is obtained for points far enough away from the free boundary. This is because, in Equation (\ref{eqn:error-greens-prop}), $G(S,S_{m+1}) = \mathcal O(S_{m+1}-S_f) = \mathcal O(h)$, and $G(S,S_{m+2}) = \mathcal O(S_{m+2}-S_f) = \mathcal O(h)$ as $h\rightarrow 0$. For a visual demonstration of property, Figure \ref{fig:bvp_greens} gives an illustration of Green's function for a hypothetical second-order differential equation on two successive grid refinements.


\subsubsection{Initial value problems}
When the time variable is included, the analysis becomes more involved. However, the conclusions are similar to the ones for boundary value problems. In this section, we study the single step error behavior when we solve \autoref{eqn:IVP-penalty-LCP} using fourth-order finite difference discretization and BDF4 time-stepping. While the stability analysis is important, we observed empirically that our time-stepping scheme is generally stable in practice. Hence, we leave the stability analysis for future research, and only focus on single step error behavior.

Consider the original LCP given by (\ref{eqn:IVP-LCP}), and the corresponding penalty equation (\ref{eqn:IVP-penalty-LCP}).
As in Proposition \ref{prop:fourth-order-fd-error} we have relation (\ref{eqn:dt-L-exact-sub}), it is easy to see that, in the case of time-dependent problems, we have
\begin{equation}\label{eqn:exact-LCP-fd}
    \partial_t\hat{\mathbf V} = \mathbf L\hat{\mathbf V} + \mathbf b + \bm{\upgamma} + \bm{\uptheta},\quad \text{ with}\; [\bm{\upgamma}(t)]_j \equiv \left((\partial_t-\mathcal L)\hat{V}(t,S_j)-g(t,S_j)\right)\mathbbm 1_{1\leq j \leq m(t)},
\end{equation}
where $m(t)$ is the node index at time $t$ such that $S_{m(t)} \leq S_f(t) < S_{m(t)+1}$.

Consider the BDF4 discretization applied to Equations
(\ref{eqn:exact-LCP-fd}) starting at the fourth time step. We have,
for the exact LCP solution,
\begin{equation}\label{eqn:BDF4-exact-LCP}
    \mathbf A\hat{\mathbf V}^{n+4} = 4\hat{\mathbf V}^{n+3} - 3\hat{\mathbf V}^{n+2} + \frac{4}{3}\hat{\mathbf V}^{n+1} - \frac{1}{4}\hat{\mathbf V}^{n} + k\mathbf b^{n+4}+ k(\bm{\upgamma} + \bm{\uptheta} + \bm{\upbeta}),
\end{equation}
where $\bm{\upbeta}$ is the truncation error of the BDF4 time-stepping scheme applied to (\ref{eqn:exact-LCP-fd}). Note that the fully discrete system that we are actually solving is Equation (\ref{eqn:BDF4-discrete}).
The following proposition gives the relationship between the solution $\hat{V}$ of the exact LCP and the solution $\tilde{\mathbf V}^{n+4}$ of the fully discrete system \autoref{eqn:BDF4-discrete}. To simplify the notation, in the proposition as well as the theorem following, we drop the superscript $n+4$ from the computed solution $\tilde{\mathbf V}^{n+4}$. Note that we have also dropped the superscript from $\bm{\upgamma}, \bm{\uptheta}$ and $\bm{\upbeta}$ for simplicity.

\begin{proposition}\label{prop:ivp-error-prop}
Let $\hat{V}$ be the solution to the exact LCP (\ref{eqn:IVP-LCP}), and $V$ be the solution to the continuous penalty equation (\ref{eqn:IVP-penalty-LCP}). Assume that
the penalty terms in the fully discrete equations reflect the correct behavior, i.e., $\mathcal I_{\tilde{\mathbf V}} = \mathcal I_{\hat{\mathbf V}}$ at each time step. Then, we have
\begin{equation}\label{eqn:BDF4-fully-discrete-erro-sysr}
    \left(\frac{25}{12k}\mathbf I - \mathbf L + \rho\mathcal I_{\hat{\mathbf V}} \right)\mathbf e^{n+4} = \frac{1}{k}\left( 4\mathbf e^{n+3} - 3\mathbf e^{n+2} + \frac{4}{3}\mathbf e^{n+1} - \frac{1}{4}\mathbf e^{n}\right) + (\bm{\upgamma}+\bm{\uptheta} + \bm{\upbeta}),
\end{equation}
where $\mathbf e = \hat{\mathbf V}-\tilde{\mathbf V}$ is error in the solution.
\end{proposition}
\begin{proof}
First, notice that $\rho\mathcal I_{\tilde{\mathbf V}}(\mathbf V^* - \tilde{\mathbf V}) = \rho\mathcal I_{\hat{\mathbf V}}(\hat{\mathbf V} - \tilde{\mathbf V})$, since
\begin{equation*}
    \rho\mathcal I_{\tilde{\mathbf V}}(\mathbf V^* - \tilde{\mathbf V}) = \rho\mathcal I_{\tilde{\mathbf V}}(\mathbf V^* - \hat{\mathbf V} + \hat{\mathbf V} - \tilde{\mathbf V}) = \rho\mathcal I_{\hat{\mathbf V}}(\mathbf V^* - \hat{\mathbf V}) + \rho\mathcal I_{\tilde{\mathbf V}}(\hat{\mathbf V} - \tilde{\mathbf V})= \rho\mathcal I_{\tilde{\mathbf V}}(\hat{\mathbf V} - \tilde{\mathbf V}) = \rho\mathcal I_{\hat{\mathbf V}}(\hat{\mathbf V} - \tilde{\mathbf V}).
\end{equation*}
Using this identity, we subtract Equation (\ref{eqn:BDF4-discrete}) from (\ref{eqn:BDF4-exact-LCP}) and rearrange to get Equation (\ref{eqn:BDF4-fully-discrete-erro-sysr}).
\end{proof}

Corresponding to Proposition \ref{prop:fourth-order-fd-error} in the previous section which describes the error equation for time-independent free boundary problems, Proposition \ref{prop:ivp-error-prop} gives an expression of the error evolution for solving time-dependent free boundary problems using a fourth-order finite difference scheme and BDF4. For convenience of discussion, define
\begin{equation}\label{eqn:Lk_rk}
    \mathbf L_k \equiv \frac{25}{12k}\mathbf I - \mathbf L, \text{ and } \mathbf r_k^{n+3} \equiv \frac{1}{k}\left( 4\mathbf e^{n+3} - 3\mathbf e^{n+2} + \frac{4}{3}\mathbf e^{n+1} - \frac{1}{4}\mathbf e^{n}\right) + (\bm{\upgamma}+\bm{\uptheta} + \bm{\upbeta}).
\end{equation}
In addition, assume that, at the $(n+4)$-th time step, the free boundary is located in between $S_m$ and $S_{m+1}$, i.e., $S_{m} \leq S_f(t_{n+4}) < S_{m+1}$, on the space grid $\{S_0, S_1, \ldots, S_m, \ldots S_M\}$. Similar to the discussion for boundary value problems, we divide the matrix $\mathbf L_k$ into four block submatrices
\begin{equation*}
    \mathbf L_k =
    \begin{bmatrix}
    [\mathbf L_k]_{11} & [\mathbf L_k]_{12}\\
    [\mathbf L_k]_{21} & [\mathbf L_k]_{22}
    \end{bmatrix}
\end{equation*}
where the block matrices $[\mathbf L_k]_{11},\; [\mathbf L_k]_{12},\;[\mathbf L_k]_{21},\; [\mathbf L_k]_{22}$ are of sizes $m\times m$, $m\times (M-m)$, $(M-m)\times m$ and $(M-m)\times(M-m)$, respectively. Corresponding to Theorem \ref{thm:penalty-error-equation} for time-independent problems, the error decomposition for time-dependent problems is given by Theorem \ref{thm:fully-discrete-error-equation}.

\begin{theorem}\label{thm:fully-discrete-error-equation}
Assume that $[\mathbf L_k]_{11}$ and $[\mathbf L_k]_{22}$ are nonsingular, and
$\rho$ is a positive number such that
$\rho \gg \max_{ij}\{|[\mathbf L_k]_{i,j}|\}$. Assume also that 
      $\max_{j}  \{|[[\mathbf L_k]_{22}^{-1}]_{1,j}|\} = \mathcal O(h^2)$,
      $\max_{j}  \{|[[\mathbf L_k]_{22}^{-1}]_{2,j}|\} = \mathcal O(h^2)$,
      $\max_{i}  \{|[[\mathbf L_k]_{22}^{-1}]_{i,1}|\} = \mathcal O(h^2)$,
and   $\max_{i}  \{|[[\mathbf L_k]_{22}^{-1}]_{i,2}|\} = \mathcal O(h^2)$. Further, assume that $\mathcal I_{\tilde{\mathbf V}} = \mathcal I_{\hat{\mathbf V}}$ at each time step. When using BDF4 time-stepping and the fourth-order finite difference scheme to solve the penalized PDE \autoref{eqn:IVP-penalty-LCP}, the solution error at the $(n+4)$-th time step satisfies
\begin{equation*}
    \mathbf e^{n+4} \approx
    \begin{bmatrix}
    \mathbf O\\
    \left([\mathbf L_{k}]_{22}\right)^{-1}[\mathbf r_{k}^{n+3}]_{m+1:M}
    \end{bmatrix}.
\end{equation*}
\end{theorem}
\begin{remark1}
We denote $([\mathbf L_k]_{22})^{-1}$ by $[\mathbf L_k]_{22}^{-1}$ for notation simplicity.
\end{remark1}
\begin{proof}
From Proposition \ref{prop:ivp-error-prop}, we know that the error is the solution to
\begin{equation*}
    (\mathbf L_k + \rho\mathcal I_{\tilde{\mathbf V}})\mathbf e^{n+4} = \mathbf r_k^{n+3}.
\end{equation*}
Therefore, we get
\begin{equation*}
    \mathbf e^{n+4} = (\mathbf L_k + \rho\mathcal I_{\tilde{\mathbf V}})^{-1}\mathbf r_k^{n+3}.
\end{equation*}
Then by making use of Proposition \ref{prop:penalized-matrix-inverse} applied to $\mathbf L_k$, the theorem is proved.
\end{proof}

\begin{remark1}
Note that similar to \autoref{prop:penalized-matrix-inverse}, the assumptions on $[\mathbf L_k]_{22}^{-1}$ are typically true in practice.
\end{remark1}

Theorem \ref{thm:fully-discrete-error-equation} shows that the errors in the approximate solutions of moving boundary problems behave in a similar way to the solution of free boundary problems. The solution on the penalty region is computed exactly within a tolerance, while on the PDE region, in addition to the truncation errors, the solution errors from previous time steps also contribute to the solution error at the current time step. The error propagation is governed by $[\mathbf L_k]_{22}^{-1}$, which depends on both the time stepping and space discretization schemes. Similar to the discussion of boundary value problems, it can be treated as the discrete analogue of the Green's function to the continuous operator on the PDE region.

\begin{figure}[h]
\begin{tabular}{c c}
    \includegraphics[width=0.45\textwidth]{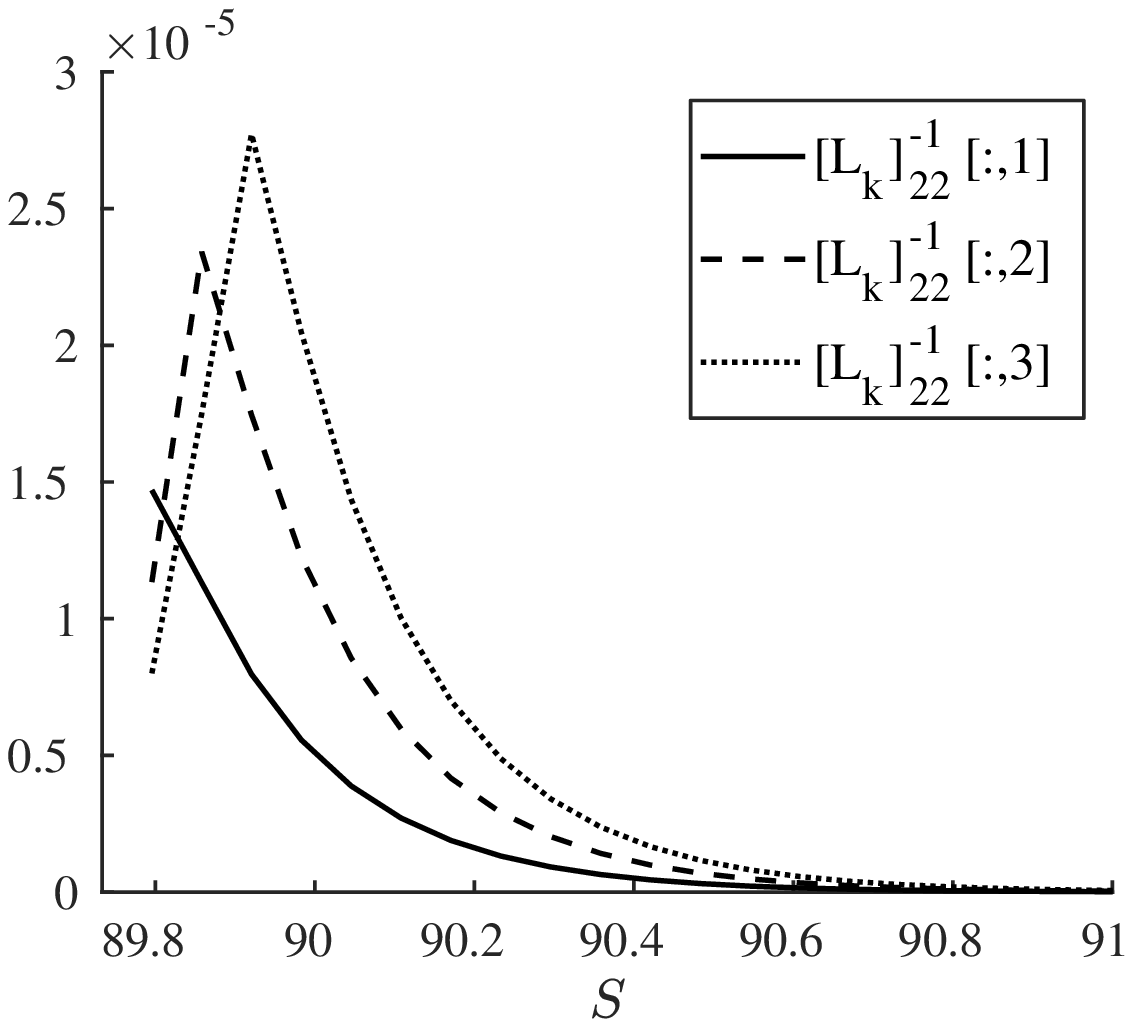} &  \includegraphics[width=0.45\textwidth]{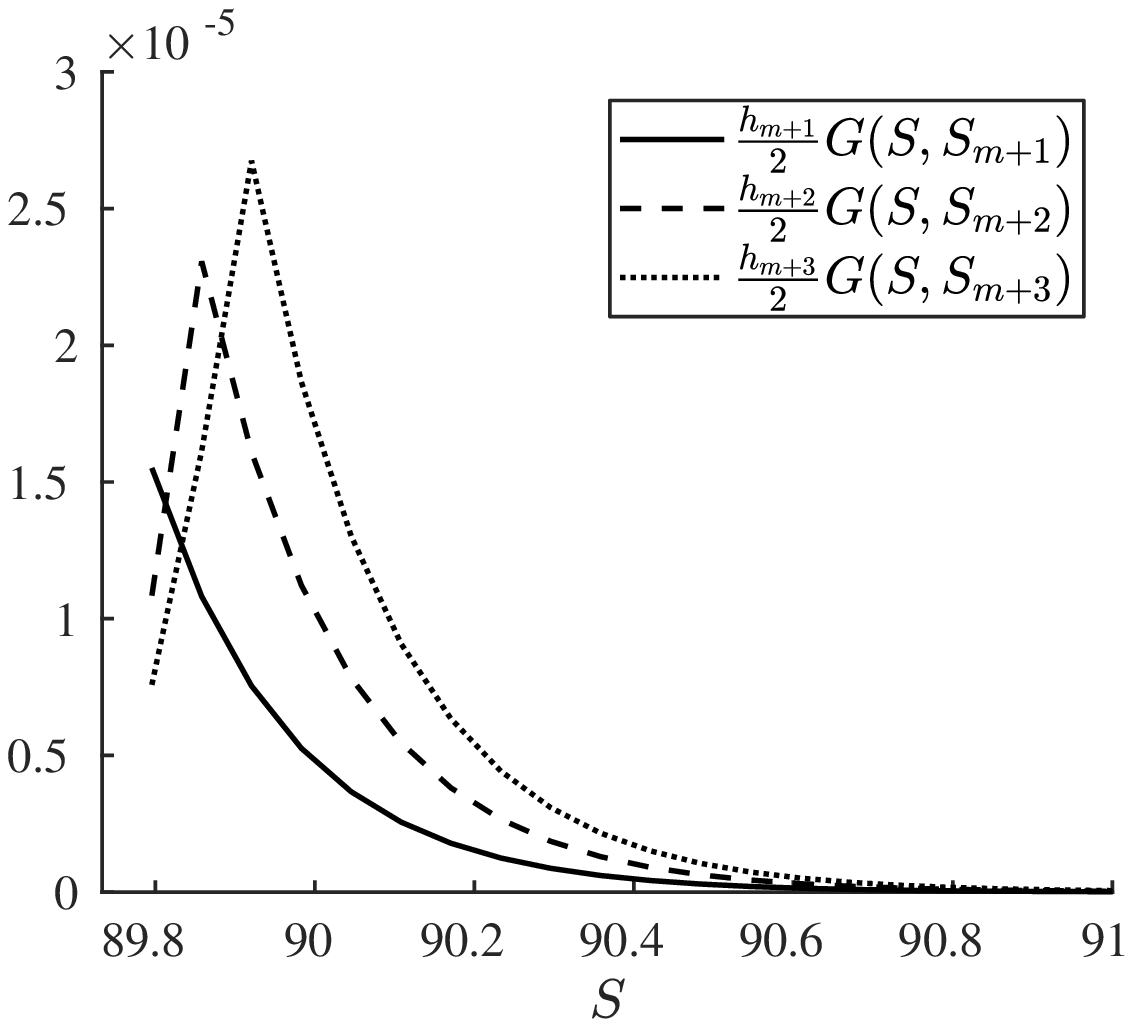}\\
    (a) The first 3 columns of $[\mathbf L_k]_{22}^{-1}$ & (b)
    The continuous Green's function $\frac{h_j}{2}G(S,S_j)$
\end{tabular}
    \centering
    \caption{(a) The first three columns of $[\mathbf L_k]_{22}^{-1}$ on an example nonuniform grid; (b) The scaled continuous Green's function $\frac{h_j}{2}G(S,S_j)$ for the operator $\mathcal L_{BS}$ at $S_{m+1}$, $S_{m+2}$ and $S_{m+3}$. The free boundary location is $S_f = 89.748$. Note that the zero value on the left of the free boundary is not included.}
    \label{fig:Lk22inv}
\end{figure}

In the following, we consider a concrete example of the Black-Scholes operator, i.e., let $\mathcal L = \mathcal L_{BS}$. Instead of studying $\mathbf L_k$, which is hard to analyze, we investigate the Green's function for the continuous operator $\mathcal L_k = \frac{25}{12k} - \mathcal L$ on the PDE region for a fixed time step size $k$. In Figure \ref{fig:Lk22inv}, we show the comparison of the graphs of $G(\mathbf S_2, S_j)$ on the PDE region and $[\mathbf L_k]_{22}^{-1}$ for the first three Green's functions on an example grid, where the free boundary is located at $S_f = 89.748$. Again, we see that they have the same shape with similar magnitudes. By performing the usual variable transformation $S = Ke^x$ to $\mathcal L_k$,
we can get a transformed operator $\mathcal L_{k,x} = -\partial_{xx} - (\kappa-1)\partial_x + (\kappa + \frac{25}{6\kappa\sigma^2})$,
whose Green's function is given by Equation \autoref{eqn:general-greens-function} in Proposition \ref{prop:greens}, with $\xi_1 =\frac{-(\kappa-1)+\sqrt{(\kappa+1)^2+4\lambda}}{2}$, $\xi_2 = \frac{-(\kappa-1)-\sqrt{(\kappa+1)^2+4\lambda}}{2}$, where $\lambda = \frac{25}{6\kappa\sigma^2}$, and $\kappa = \frac{2r}{\sigma^2}$. Let $x_{m+1} = \log(S_{m+1}/K)$ and $x_f = \log(S_f/K)$. When $x_{m+1}-x_f$ is small enough, we have that $G(x,x_{m+1}) = \mathcal O(x_{m+1}-x_f)$ by Proposition \ref{prop:greens}. Since
\begin{equation*}
    x_{m+1}-x_f = \log(S_{m+1}/K) - \log(S_f/K) = \log\left(1+\frac{S_{m+1}-S_f}{S_f}\right) \approx \frac{S_{m+1}-S_f}{S_f} = \mathcal O(h),
\end{equation*}
when $S_{m+1}\rightarrow S_f > 0$, the Green's function $G(S,S_{m+1}) = \mathcal O(h)$.
Therefore, using the same argument as in Section \ref{subsec:bvp-conv}, we claim that the error in the solution of Equation \autoref{eqn:BDF4-discrete} is of order $\mathcal O(h^2)$.

\subsection{Grid crossing}\label{sec:time-discretization}
With BDF4 time-stepping, the time derivative of the solution at some point is computed by a linear combination of the solutions at the four points directly prior to the current point. However, one of more of the prior points may not lie on the same side of the moving boundary. 

\begin{figure}[h]
    \centering
    \includegraphics[width=0.65\textwidth]{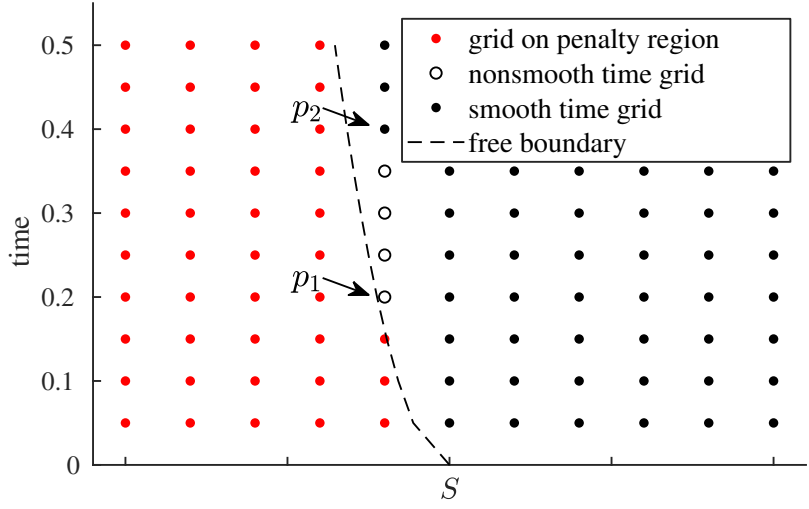}
    \caption{{\bf An example layout of grid points in the time and space domain.}
    The dashed line is the free boundary. The red points to the left of the free boundary are on the penalty region, and
    the points to the right of the free boundary are on the PDE region.
    Unlike at the solid black points, BDF4 has degenerated accuracy at the hollow black points because it involves solution points that lie on different sides of the free boundary.}
    \label{fig:grid_labels}
\end{figure}

To see this, we can look at an example grid shown in Figure \ref{fig:grid_labels}. The black hollow points, for example $p_1$, are problematic points for the BDF4 time-stepping scheme, because their computation depends on one or more points on the other side of the moving boundary. Therefore, the BDF4 scheme at the black hollow points may exhibit degenerated accuracy. 

On the other hand, for the black solid points, such as $p_2$, the time derivative of the solution with the BDF4 scheme uses only points on the same side of the moving boundary. Hence, the BDF4 scheme at the black solid points is fourth-order accurate. In general, we can see that the number of problematic points depends on how quickly the moving boundary moves relative to the grid spacing. On a fixed grid, a slow-moving free boundary will have a smaller number of black hollow points.

 At minimum, we require the time discretization to be accurate enough so as not to affect the convergence order in space. Despite the complicated behavior of the time derivative discretization, we do not explicitly deal with the loss of accuracy in the time derivative in our algorithm. Instead, we apply a time variable transformation to change the shape of the free boundary and try to reduce the number of problematic points in the time stepping. This time transformation together with appropriate space stretching turn out to be good enough to maintain high-order accuracy.

\subsection{Extrapolating the numerical solution}\label{sec:extrap}
In this section, we discuss how to approximate the derivative jumps and the free boundary location using a given numerical solution, and analyze some related technical details.

\subsubsection{Approximation of the solution derivatives and associated derivative jumps}\label{sec:deriv_approx}

In this subsection, we discuss how to approximate the derivative jumps
$\Delta V_{S_f}''   \equiv V_{S_f,-}''   - V_{S_f,+}''$,
$\Delta V_{S_f}'''  \equiv V_{S_f,-}'''  - V_{S_f,+}'''$, and
$\Delta V_{S_f}'''' \equiv V_{S_f,-}'''' - V_{S_f,+}''''$,
where the ``$-$'' and the ``$+$'' in the subscripts denote
values in the penalty and the PDE regions, respectively.

In the penalty region, since we are given the obstacle function $V^*$,
and we know that $V = V^*$,
we can evaluate $V_{S_f,-}''$, $V_{S_f,-}'''$ and $V_{S_f,-}''''$ exactly.

\begin{figure}
    \centering
    \includegraphics[width=0.55\textwidth]{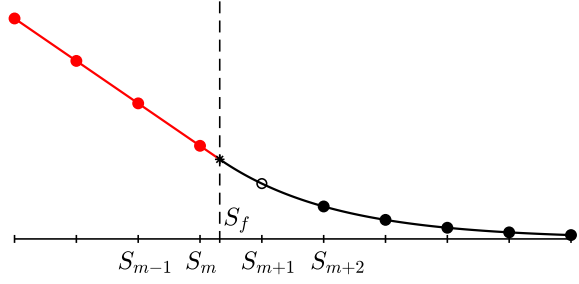}
    \caption{An example grid on which a free boundary problem is defined. Point $S_{f}$ is the free boundary location.}
    \label{fig:extrap-grid}
\end{figure}

We now turn to the PDE region and discuss the approximation of
$V_{S_f,+}''$, $V_{S_f,+}'''$ and $V_{S_f,+}''''$.
Recall that $S_m\leq S_f < S_{m+1}$, as shown on an example uniform grid in Figure \ref{fig:extrap-grid}.
From relation \autoref{eqn:error-greens-prop},
we know that the numerical solution of PDE \autoref{eqn:BVP-penalty-LCP} by a standard fourth-order solve is exact within a tolerance $\mathcal O(1/\rho)$ at points $S_j \leq S_{m}$ and second-order accurate at the points $S_j \geq S_{m+2}$.
More precisely, the error is of the form
\begin{equation}\label{eqn:error-components}
\begin{aligned}
    \mathbf e &\approx
    \mathcal O(h)\begin{bmatrix}\mathbf 0\\G(\mathbf S_2, S_{m+1})\end{bmatrix} + \mathcal O(h)\begin{bmatrix}\mathbf 0\\G(\mathbf S_2, S_{m+2})\end{bmatrix} + \sum_{j=m+3}^M\mathcal O(h^5)\begin{bmatrix}\mathbf 0\\G(\mathbf S_2, S_{j})\end{bmatrix}.
\end{aligned}
\end{equation}
Since the Green's functions $G(S, S_j)$ are piecewise smooth with first-derivative jumps at points $S_j$, and the summation terms in \autoref{eqn:error-components} are negligible compared to the $\mathcal O(h)$ terms, the solution error $\mathbf e$ is of order $\mathcal O(h^2)$,
\begin{equation*}
    e_j \approx \mathcal O(h)G(S_j, S_{m+1}) + \mathcal O(h)G(S_j, S_{m+2}),
\end{equation*}
and it is smooth for $S_j \geq S_{m+2}$, which means the numerical solution is also smooth. This is a key to the success of our method to approximate derivatives in the PDE region without of loss of accuracy. In order to maintain the same order of accuracy as the numerical solutions, when applying finite differences to approximate the derivatives, we should only use points on the right of $S_{m+2}$ (inclusive).
Therefore, by having $O(h^2)$ accurate values $\tilde{V}_{m+2}, \tilde{V}_{m+3}, \tilde{V}_{m+4}$, we construct $O(h^2)$ accurate second derivative values $\tilde{V}''_{m+2}, \tilde{V}''_{m+3}, \tilde{V}''_{m+4}$. To approximate the second derivative of the solution at the free boundary $S_f < S_{m+1}$ to $\mathcal O(h^2)$ accuracy,
we extrapolate the solution derivative at $S_f$ using the computed $\tilde{V}''_{m+2}, \tilde{V}''_{m+3}, \tilde{V}''_{m+4}$ by
\begin{equation*}
    \tilde{V}''(S) = \sum_{i=2}^4\left(\tilde{V}''_{m+i}\prod_{j=2,j\neq i}^4\frac{S-S_{m+j}}{S_{m+i}-S_{m+j}}\right),
\end{equation*}
i.e., by using a quadratic Lagrange polynomial. The obtained approximate derivative at $S_f$ is of $\mathcal O(h^2)$ accuracy. Higher-order derivatives are computed in the same way. Hence, from the analysis in Section \ref{sec:fd-w-approx-correction}, we see that the correction terms computed using the approximate derivative jumps will be of $\mathcal O(h^2)$ accurate, which is more than the required $\mathcal O(h)$ accuracy to increase the order of the corrected finite differences. Similarly, we use $\tilde{V}_{m+2}, \tilde{V}_{m+3}, \tilde{V}_{m+4}, \tilde{V}_{m+5}$ and a cubic polynomial for the third-order approximation, and $\tilde{V}_{m+2}, \tilde{V}_{m+3}, \tilde{V}_{m+4}, \tilde{V}_{m+5}, \tilde{V}_{m+6}$ and a quartic polynomial for the fourth-order approximation of the solution derivatives at the free boundary. We stress here that the choice of the interpolation points starting from $S_{m+2}$ (skipping $S_{m+1}$) to maintain the convergence order is guided by our analysis of the error behavior in the solution given by Equation \autoref{eqn:error-components}. Similar interpolation scheme can be found in \cite{linnick2005high}, where the authors are only able to justify their results empirically from numerical experiments. Note that we could have used one degree less Lagrange interpolation in each case, but our choice is dictated by the fact that we want the extrapolation error to be of even lower order than the solution error.


\subsubsection{Approximation of the free boundary location}\label{sec:fb-approximation}
To approximate the free boundary location, we apply the smooth pasting condition of the derivative at the free boundary in Equation (\ref{eqn:value-deriv-matching}), which we repeat here for convenience:
\begin{equation}\label{eqn:smooth-pasting}
    \frac{\partial V}{\partial S}(t, S_f(t)) = \frac{\partial V^*}{\partial S}(S_f(t)).
\end{equation}
Assume that we have already calculated the approximate derivatives $\tilde{V}'_{m+2}$, $\tilde{V}'_{m+3}$, $\tilde{V}'_{m+4}$, $\tilde{V}'_{m+5}$ and $\tilde{V}'_{m+6}$ at $S_{m+2}$, $S_{m+3}$, $S_{m+4}$, $S_{m+5}$, $S_{m+6}$, respectively, with certain accuracy, using solution values starting from $S_{m+2}$ and on. We then apply the quartic Lagrange polynomial to fit the derivative by
\begin{equation*}
    \tilde{V}'(S) = \sum_{i=2}^6\left(\tilde{V}'_{m+i}\prod_{j=2,j\neq i}^6\frac{S-S_{m+j}}{S_{m+i}-S_{m+j}}\right).
\end{equation*}
Then, the approximate free boundary is obtained by Newton's root finding algorithm such that
\begin{equation*}
    \tilde{V}'(S) - \frac{\partial V^*}{\partial S}(S) \approx 0.
\end{equation*}
The approximate free boundary obtained in this way is expected to be of the same order as the numerical solution.
Hence, from the analysis in Section \ref{sec:fd-w-approx-correction}, by using $\mathcal O(h^2)$ accurate solution, the approximate free boundary is $\mathcal O(h^2)$, and the correction terms computed are $\mathcal O(h)$, as is required to increase the order of the corrected finite differences.
Similarly, using $\mathcal O(h^3)$ accurate solution, the correction terms computed are $\mathcal O(h^2)$, and so on. Note that we decided to use the smoothing pasting condition (\ref{eqn:smooth-pasting}) to locate the free boundary instead of the value matching condition $V(t,S_f(t)) = V^*(S_f(t))$. The reason for this is that the value matching equation has a zero derivative at the root. Therefore, Newton's root-finding method is slow if value matching is used.

\section{Algorithm}\label{sec:algorithm}

\subsection{A fourth-order deferred correction algorithm for solving free boundary problems}
We are now ready to present a fourth-order deferred correction finite difference algorithm for solving free and moving boundary problems. To start, we first present the algorithm for solving free boundary problems where no time variable is involved. 
Recall that the penalized equation for solving free boundary problems is given by Equation (\ref{eqn:BVP-penalty-LCP}),
and the respective discrete equations are given by the nonlinear system \autoref{eqn:BVP-discrete}, solved by a generalized Newton's iteration as described in  \cite{forsyth2002quadratic}, which is also referred to as discrete penalty iteration.

The main idea of our algorithm is to use a deferred correction technique to eliminate the lower-order errors in the finite difference approximation introduced by piecewise smoothness in the solution.
We illustrate our correction scheme by considering only the leading order terms of the corrections in Equations (\ref{eqn:correction-1})--(\ref{eqn:correction-4}), that is, those associated with the jump in the second derivative. The other terms are corrected in the same manner and we omit the discussion.

Recall that the fourth-order FD discretization of $\frac{\partial^2V}{\partial S^2}$ and $\frac{\partial V}{\partial S}$ is $\bar{\mathbf L}_2\tilde{\mathbf V}_{\text{aug}}$ and $\bar{\mathbf L}_1\tilde{\mathbf V}_{\text{aug}}$, respectively, at the interior nodes of a grid $S_0 < S_1 < \ldots < S_{M+1}$, where $\bar{\mathbf L}_2$ and $\bar{\mathbf L}_1$ are $M\times (M+2)$ second- and first-derivative, respectively, FD coefficient matrices, and $\tilde{\mathbf V}_{\text{aug}} = [\tilde{V}_0, \tilde{V}_1, \ldots, \tilde{V}_{M+1}]^T$, as defined in Section \ref{sec:discretization}. Suppose that 
$S_m \leq S_f < S_{m+1}$, as shown in Figure \ref{fig:extrap-grid}. The second derivative jump at the free boundary is pre-computed to be $\Delta V_{S_f}''$ (either approximate or exact). From Theorem \ref{thm:fourth-order-correction}, making $\delta = S_{m+1}-S_f$, and picking appropriate FD coefficients,
we see that the correction terms corresponding to the second derivative jumps
at nodes $S_{m-1},\; S_m,\; S_{m+1}$ and $S_{m+2}$,
are computed by
\begin{align}
    &C''_{1,0} = \frac{(S_{m+1}-S_f)^2}{2}\bar{\mathbf L}_2(m-1,m+2)\Delta V_{S_f}'',\\
    &C''_{2,0} = \left(\frac{(S_{m+1}-S_f)^2}{2}\bar{\mathbf L}_2(m,m+2) + \frac{(S_{m+2}-S_f)^2}{2}\bar{\mathbf L}_2(m,m+3)\right)\Delta V_{S_f}'',\\
    &C''_{3,0} = -\left(\frac{(S_f-S_{m-1})^2}{2}\bar{\mathbf L}_2(m+1,m) + \frac{(S_f-S_m)^2}{2}\mathbf L_2(m+1,m+1)\right)\Delta V_{S_f}'',\\
    &C''_{4,0} = -\frac{(S_f-S_m)^2}{2}\bar{\mathbf L}_2(m+2,m+1)\Delta V_{S_f}'',
\end{align}
where $\bar{\mathbf L}_2(i,j)$ denotes the $(i,j)$ entry of the coefficient matrix $\bar{\mathbf L}_2$.
Note that $C''_{j,0}$, $j = 1, \ldots, 4$, correspond
to $\mathcal O(1)$ error terms in
(\ref{eqn:correction-1})-(\ref{eqn:correction-4}).
The correction terms corresponding to the third and fourth derivative jumps
are computed similarly, giving rise to
$C''_{j,1}$, $j = 1, \ldots, 4$ (corresponding to $\mathcal O(h)$
error terms in (\ref{eqn:correction-1})-(\ref{eqn:correction-4})) and
$C''_{j,2}$, $j = 1, \ldots, 4$ (corresponding to $\mathcal O(h^2)$
error terms in (\ref{eqn:correction-1})-(\ref{eqn:correction-4})),
respectively.
Then, the total correction terms for the FD approximation
of the second derivative
at nodes $S_{m-1},\; S_m,\; S_{m+1}$ and $S_{m+2}$ are
\begin{align}
    &C''_1 = C''_{1,0} + C''_{1,1} + C''_{1,2},\\
    &C''_2 = C''_{2,0} + C''_{2,1} + C''_{2,2},\\
    &C''_3 = C''_{3,0} + C''_{3,1} + C''_{3,2},\\
    &C''_4 = C''_{4,0} + C''_{4,1} + C''_{4,2}.
\end{align}
By replacing the $\bar{\mathbf L}_2$ entries with the corresponding $\bar{\mathbf L}_1$ entries, we can similarly calculate the correction terms $C'_j$ to the first derivative approximations, for $j = 1, 2, 3, 4$, at nodes $S_{m-1},\; S_m,\; S_{m+1}$ and $S_{m+2}$, respectively.

With these correction entries in hand, instead of solving Equation (\ref{eqn:BVP-discrete}), we solve a modified system
\begin{equation}\label{eqn:BVP-discrete-corrected}
    \mathbf L\tilde{\mathbf V}+\mathbf b + \rho\mathcal I_{\tilde{\mathbf V}}(\mathbf V^*-\tilde{\mathbf V})  + \mathbf a_2+ \mathbf a_1 = \mathbf 0,
\end{equation}
with correction terms $\mathbf a_2$ and $\mathbf a_1$, where
\begin{align}
    &\mathbf a_2 = [0,\; \ldots,\; 0,\; p(S_{m-1})C''_1,\; p(S_{m})C''_2,\; p(S_{m+1})C''_3,\; p(S_{m+2})C''_4,\; 0,\; \ldots,\; 0]^T,\\
    &\mathbf a_1 = [0,\; \ldots,\; 0,\; w(S_{m-1})C'_1,\; w(S_{m})C'_2,\; w(S_{m+1})C'_3,\; w(S_{m+2})C'_4,\; 0,\; \ldots,\; 0]^T.
\end{align}
Ideally, we would know the exact derivative jumps and apply them to correct the FDs when discretizing the PDE. However, these jumps are not known a priori in the setting of this work. Therefore, we make use of the approximate solution derivatives and free boundary location that we get, as described in Section \ref{sec:extrap}, to compute the approximate correction terms.
Applying these correction terms increases the order of the finite difference approximation.
As a result, it also increases the order of the new approximate solution when we solve the discrete system again with the corrected FDs. We give the details in Algorithm \ref{alg:fb}.\\

\begin{algorithm}[H]
\caption{A fourth-order FD algorithm for solving free boundary problems}
\label{alg:fb}
\begin{algorithmic}[1]
\STATE
\textbf{Phase 1}: \\
Solve (\ref{eqn:BVP-discrete}) to obtain an approximate solution $\tilde{\mathbf V}^{(0)}$ of $\mathcal O(h^2)$ using penalty iteration.\\
Find $S_m$ and $S_{m+1}$ as in Proposition \ref{prop:fourth-order-fd-error}.\\
Compute the approximate free boundary $\tilde{S}_f^{(0)}$ to $\mathcal O(h^2)$ accuracy,  using Newton's method with initial guess $(S_m+S_{m+1})/2$, as in Section \ref{sec:fb-approximation}. \\
Approximate $V''(S_f)$ at $\tilde{S}_f^{(0)}$ to obtain ${\tilde{V}^{''(0)}}_{S_f}=V''(S_f)+O(h)$, as in Section \ref{sec:deriv_approx}.\\
\STATE
\textbf{Phase 2}: \\
Compute the FD corrections $\mathbf a_1$, $\mathbf a_2$ using $\tilde{S}_f^{(0)}$ and ${V^*_{S_f}}''-{\tilde{V}^{''(0)}}_{S_f}$. \\
Solve (\ref{eqn:BVP-discrete-corrected}) to obtain an approximate solution $\tilde{\mathbf V}^{(1)}$ of $\mathcal O(h^3)$ using penalty iteration with initial guess $\tilde{\mathbf V}^{(0)}$.\\
Compute the approximate free boundary $\tilde{S}_f^{(1)}$ to $\mathcal O(h^3)$ accuracy, using Newton's method with initial guess $\tilde{S}_f^{(0)}$. \\
Approximate $V''(S_f),\; V'''(S_f)$ at $\tilde{S}_f^{(1)}$ to obtain $\tilde{V}^{''(1)}_{S_f} = V''(S_f) + \mathcal O(h^2)$, $\tilde{V}^{'''(1)}_{S_f} = V'''(S_f) + \mathcal O(h)$.\\
\STATE
\textbf{Phase 3}: \\
Compute the FD corrections $\mathbf a_1$, $\mathbf a_2$ using $\tilde{S}_f^{(1)}$ and ${V^*_{S_f}}''-{\tilde{V}^{''(1)}}_{S_f}, \; {V^*_{S_f}}'''-{\tilde{V}^{'''(1)}}_{S_f}$.\\
Solve (\ref{eqn:BVP-discrete-corrected}) to obtain an approximate solution $\tilde{\mathbf V}^{(2)}$ of $\mathcal O(h^4)$ using penalty iteration with initial guess $\tilde{\mathbf V}^{(1)}.$\\
Compute the approximate free boundary $\tilde{S}_f^{(2)}$ to $\mathcal O(h^4)$ accuracy, using Newton's method with initial guess $\tilde{S}_f^{(1)}$. \\
Approximate $V''(S_f),\; V'''(S_f),\; V''''(S_f)$ at $\tilde{S}_f^{(2)}$ to obtain $\tilde{V}^{''(2)}_{S_f} = V''(S_f) + \mathcal O(h^3)$, $\tilde{V}^{'''(2)}_{S_f} = V'''(S_f) + \mathcal O(h^2)$, $\tilde{V}^{''''(2)}_{S_f} = V''''(S_f) + \mathcal O(h)$.\\
\STATE
\textbf{Phase 4}: \\
 Compute the FD corrections $\mathbf a_1$, $\mathbf a_2$ using $\tilde{S}_f^{(2)}$ and ${V^*_{S_f}}''-{\tilde{V}^{''(2)}}_{S_f}$, ${V^*_{S_f}}'''-{\tilde{V}^{'''(2)}}_{S_f}$, ${V^*_{S_f}}''''-{\tilde{V}^{''''(2)}}_{S_f}$.\\
 Solve (\ref{eqn:BVP-discrete-corrected}) to obtain an approximate solution $\tilde{\mathbf V}^{(3)}$ using penalty iteration with initial guess $\tilde{\mathbf V}^{(2)}$.\\
 Compute the approximate free boundary $\tilde{S}_f^{(3)}$, using Newton's method with initial guess $\tilde{S}_f^{(2)}$.
\end{algorithmic}
\end{algorithm}

\begin{remark1}
In fact, the algorithm already reaches fourth-order accuracy with two corrections. However, from numerical results presented later, it turns out the third correction improves the error noticeably. Therefore, we present the algorithm with three corrections.
\end{remark1}

\subsection{A fourth-order deferred correction algorithm for solving moving boundary problems} \label{sec:alg-mb}

When solving moving boundary problems, we also need to consider time discretization. In order to show the flow of computations as the correction phases and
timesteps proceed, we introduce a double index, with
$n$ denoting the timestep and $\ell$ (in parentheses) the correction phase.
We assume that Equation (\ref{eqn:IVP-penalty-LCP}) has been discretized in time by BDF4, except for the first three time steps, and in space by standard fourth-order FDs, resulting in the nonlinear system \autoref{eqn:BDF4-discrete}.
At this point, we make a note regarding the choice of $\rho$
in the discrete problem.
As discussed in \cite{forsyth2002quadratic},
it may be appropriate to adjust the value of $\rho$
for each refinement of the grid in a way so that
the error arising from the approximation of the LCP
by the penalized nonlinear PDE
reduces at the same rate as the discretization error.
However, it is more practical to set a small enough target
relative error tolerance $tol$
in the approximation of the LCP by the penalized nonlinear PDE.
Following the same arguments as in \cite{forsyth2002quadratic},
and under similar boundedness assumptions,
we essentially scale $\rho$ as $k^{-1}$
in (\ref{eqn:BDF4-discrete}), or, equivalently,
solve, with a fixed $\rho$, the nonlinear system
\begin{equation}\label{eqn:BDF4-discrete-noscale}
    \mathbf A\tilde{\mathbf V}^{n+4,(\ell)} = \tilde{\mathbf y}^{n+4,(\ell)}
   + \rho\bm{\mathcal I}_{\tilde{\mathbf V}^{n+4,(\ell)}}
   (\mathbf V^*-\tilde{\mathbf V}^{n+4,(\ell)}),
\end{equation}
where
\begin{equation}\label{eqn:fb-discrete-mat-rhs}
    \mathbf A = \left( \frac{25}{12}\mathbf I - k\mathbf L \right), \quad
    \tilde{\mathbf y^{n+4,(\ell)}} = k\tilde{\mathbf b}^{n+4}
   + 4\tilde{\mathbf V}^{n+3,(\ell)}
   - 3\tilde{\mathbf V}^{n+2,(\ell)}
   + \frac{4}{3}\tilde{\mathbf V}^{n+1,(\ell)}
   - \frac{1}{4}\tilde{\mathbf V}^{n,(\ell)}.
\end{equation}
When we apply corrections to Equation (\ref{eqn:BDF4-discrete-noscale}),
we solve a modified system
\begin{equation}\label{eqn:mv-discrete-corrected}
    \mathbf A\tilde{\mathbf V}^{n+4,(\ell)} = \tilde{\mathbf y}^{n+4,(\ell)}
   + \rho\bm{\mathcal I}_{\tilde{\mathbf V}^{n+4,(\ell)}}
   (\mathbf V^*-\tilde{\mathbf V}^{n+4,(\ell)})
   + k(\mathbf a_1 + \mathbf a_2).
\end{equation}
with correction terms $\mathbf a_1$ and $\mathbf a_2$ computed
in the same way as for free boundary problems discussed
in the previous section. 
Only slight modifications to Algorithm \ref{alg:fb} are required
to include BDF4 time-stepping and solve moving boundary problems.
Our high-order finite difference method for solving moving boundary
problems is given in Algorithm \ref{alg:mb}.

\begin{remark1}
It is important to note that,
to produce quantities for the $\ell$th phase at $t_{n+4}$,
data from the $\ell$th phase of timesteps
$t_{n+3}$, $t_{n+2}$, $t_{n+1}$ and $t_{n}$ are used.
\end{remark1}

\begin{algorithm}[H]
\caption{A fourth-order FD algorithm for solving moving boundary problems}
\label{alg:mb}
\begin{algorithmic}[1]
\FOR{\text{each time step} $t_n$}
\STATE
\textbf{Phase 1} ($\ell = 0$): \\
Solve (\ref{eqn:BDF4-discrete-noscale}) to obtain an approximate solution $\tilde{\mathbf V}^{n,(0)}$ of $\mathcal O(h^2)$ using penalty iteration with initial guess $\tilde{\mathbf V}^{n-1,(0)}$.\\ 
Find $S_{m}, S_{m+1}$ at $t_n$ as in Proposition \ref{prop:fourth-order-fd-error}.\\
Compute the approximate free boundary $\tilde{S}_f^{n,(0)}$ to $\mathcal O(h^2)$ accuracy, using Newton's method with initial guess $(S_m+S_{m+1})/2$ as in Section \ref{sec:fb-approximation}. \\
Approximate $V''(t_n,S_f)$ at $\tilde{S}_f^{n,(0)}$ to obtain ${\tilde{V}^{''(0)}}_{S_f}=V''(t_n,S_f)+O(h)$, as in Section \ref{sec:deriv_approx}.\\
\STATE
\textbf{Phase 2} ($\ell = 1$): \\
Compute the FD corrections $\mathbf a_1$, $\mathbf a_2$ using $\tilde{S}_f^{n,(0)}$ and ${V^*_{S_f}}''-{\tilde{V}^{''n,(0)}}_{S_f}$.\\
Solve (\ref{eqn:mv-discrete-corrected}) to obtain an approximate solution $\tilde{\mathbf V}^{n,(1)}$ of $\mathcal O(h^3)$ using penalty iteration with initial guess $\tilde{\mathbf V}^{n,(0)}$.\\
Compute the approximate free boundary $\tilde{S}_f^{n,(1)}$ to $\mathcal O(h^3)$ accuracy, using Newton's method with initial guess $\tilde{S}_f^{n,(0)}$. \\
Approximate $V''(t_n,S_f),\; V'''(t_n,S_f)$ at $\tilde{S}_f^{n,(1)}$ to obtain $\tilde{V}^{''n,(1)}_{S_f} = V''(t_n,S_f) + \mathcal O(h^2), \; \tilde{V}^{'''n,(1)}_{S_f} = V'''(t_n,S_f) + \mathcal O(h)$.\\
\STATE
\textbf{Phase 3} ($\ell = 2$): \\
Compute the FD corrections $\mathbf a_1$, $\mathbf a_2$ using $\tilde{S}_f^{n,(1)}$ and ${V^*_{S_f}}''-{\tilde{V}^{''n,(1)}}_{S_f}, \; {V^*_{S_f}}'''-{\tilde{V}^{'''n,(1)}}_{S_f}$.\\
Solve (\ref{eqn:mv-discrete-corrected}) to obtain an approximate solution $\tilde{\mathbf V}^{n,(2)}$ of $\mathcal O(h^4)$ using penalty iteration with initial guess $\tilde{\mathbf V}^{n,(1)}$.\\
Compute the approximate free boundary $\tilde{S}_f^{n,(2)}$ to $\mathcal O(h^4)$ accuracy, using Newton's method with initial guess $\tilde{S}_f^{n,(1)}$.\\
Approximate $V''(t_n,S_f)$, $V'''(t_n,S_f)$, $V''''(t_n,S_f)$ at $\tilde{S}_f^{n,(2)}$ to obtain $\tilde{V}^{''n,(2)}_{S_f} = V''(t_n,S_f) + \mathcal O(h^3)$, $\tilde{V}^{'''n,(2)}_{S_f} = V'''(t_n,S_f) + \mathcal O(h^2)$, $\tilde{V}^{''''n,(2)}_{S_f} = V''''(t_n,S_f) + \mathcal O(h)$.\\
\STATE
\textbf{Phase 4} ($\ell = 3$): \\
Compute the FD corrections $\mathbf a_1$, $\mathbf a_2$ using $\tilde{S}_f^{n,(2)}$ and ${V^*_{S_f}}''-{\tilde{V}^{''n,(2)}}_{S_f}$, ${V^*_{S_f}}'''-{\tilde{V}^{'''n,(2)}}_{S_f}$, ${V^*_{S_f}}''''-{\tilde{V}^{''''n,(2)}}_{S_f}$.\\
Solve (\ref{eqn:mv-discrete-corrected}) to obtain an approximate solution $\tilde{\mathbf V}^{n,(3)}$ using penalty iteration with initial guess $\tilde{\mathbf V}^{n,(2)}$.\\
Compute the approximate free boundary $\tilde{S}_f^{n,(3)}$, using Newton's method with initial guess $\tilde{S}_f^{n,(2)}$.\\
\ENDFOR
\end{algorithmic}
\end{algorithm}

\vspace{3mm}

\section{Numerical results}\label{sec:results}
In this section, we demonstrate the effectiveness of our high-order deferred correction algorithm for solving free boundary problems with several examples.
We begin with a simple one-dimensional elliptic obstacle problem in Section \ref{sec:bvp-numerical-example}. We then consider a time-dependent problem in Section \ref{sec:simple-ivp-numerical-results}. In this problem, the free boundary position $x_f(t)$ has first derivative singularity at $t = 0$, but the solution itself is smooth everywhere. For both problems, we show that the convergence rate at each solve phase in Algorithm \ref{alg:fb} and \ref{alg:mb}, respectively, is exactly as predicted.
Finally, in Section \ref{sec:american-numerical-results}, we apply our Algorithm \ref{alg:mb} with some additional considerations to the American option pricing problem. For this problem, both the free boundary and the solution itself have first  derivative singularity at $t = 0$. The numerical results show the expected convergence rate at each solve phase. 

\subsection{Solving an elliptic obstacle problem}\label{sec:bvp-numerical-example}
Consider the boundary value problem defined in the LCP form
\begin{equation}\label{eqn:bvp_ex}
\begin{aligned}
-f'' + f + 1&\geq 0,\\
f - f^* &\geq 0,\\
(-f'' + f + 1 = 0) \vee (f-f^*&=0),
\end{aligned}
\end{equation}
on the domain $x\in[-1,1]$, where $f^*(x) = x$, with boundary conditions
\begin{equation*}
    f(-1) = -1, \; f(1) = e-1.
\end{equation*}
The exact solution to this problem is the piecewise smooth function
\begin{equation*}
    f(x) = \left\{
    \begin{array}{lr}
    e^x - 1, & 0 < x \leq 1,\\
    x, & -1 \leq x \leq 0.
    \end{array}\right.
\end{equation*}
It is obvious that the solution is smooth on both $[-1, 0]$ and $[0,1]$, separately. At the point $x = 0$, the solution satisfies the value matching and smooth pasting conditions, that is,
\begin{align*}
    \lim_{x\rightarrow0^-} f(x) = \lim_{x\rightarrow0^+} f(x) = f(0) = 0, \text{ and } \lim_{x\rightarrow0^-} f'(x) = \lim_{x\rightarrow0^+} f'(x) = f'(0) = 1.
\end{align*}
However, the solution $f(x)$ has a discontinuous second derivative at $x = 0$, which means $f(x)\in C^1\backslash C^2$.
To apply Algorithm \ref{alg:fb}, we write (\ref{eqn:bvp_ex}) in the penalty form
\begin{equation}\label{eqn:bvp_ex_penalty}
    -f'' + f + 1 - \rho\max(f^*-f,0) = 0,
\end{equation}
with Dirichlet boundary conditions at $x = -1$ and 1, where $\rho$ is a penalty constant, taken to be $\rho = 1\times 10^{12}$ in the numerical experiments.

\begin{figure}[h]
    \centering
\begin{tabular}{c c}
    \includegraphics[width=0.48\textwidth]{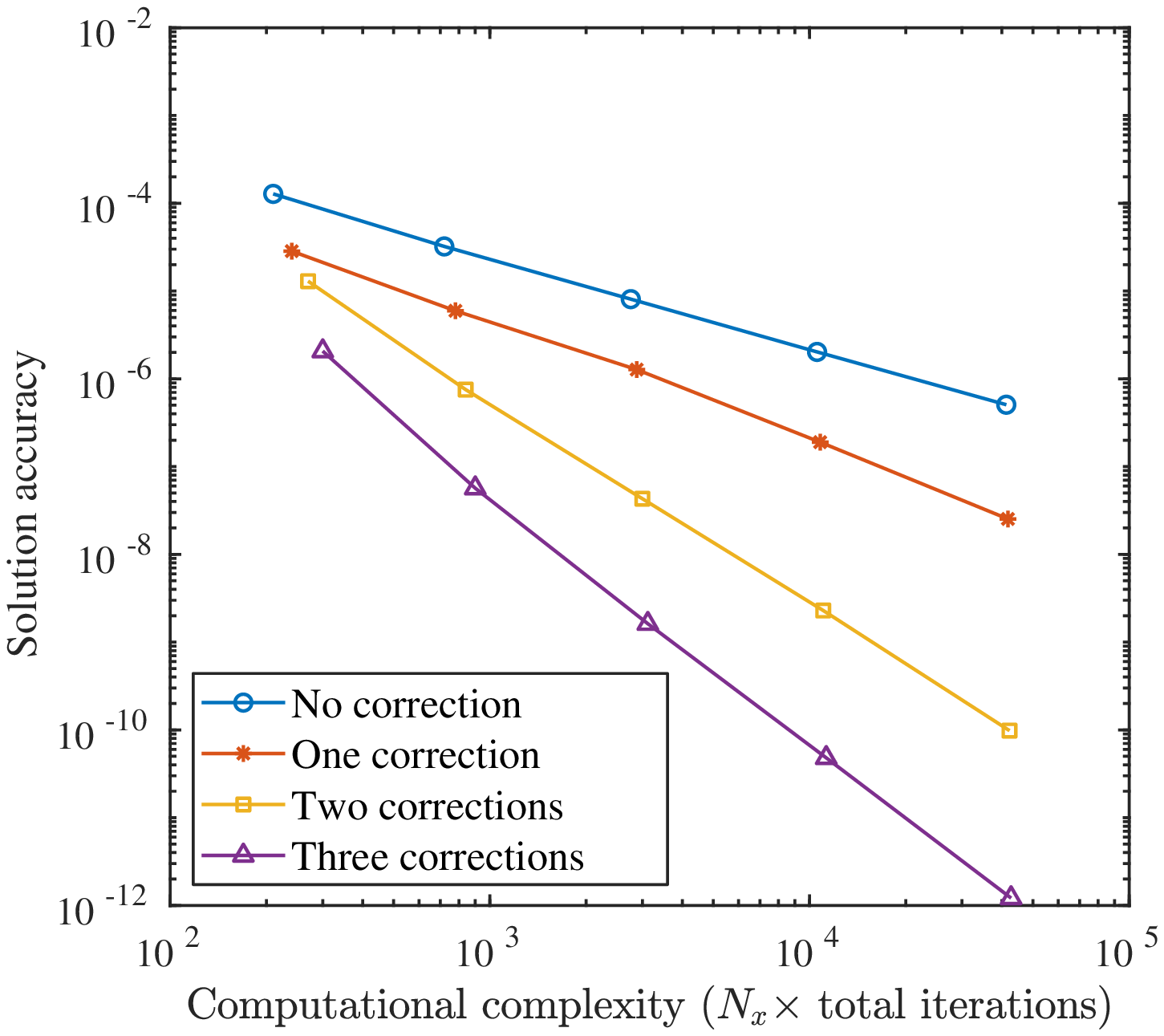} &  \includegraphics[width=0.48\textwidth]{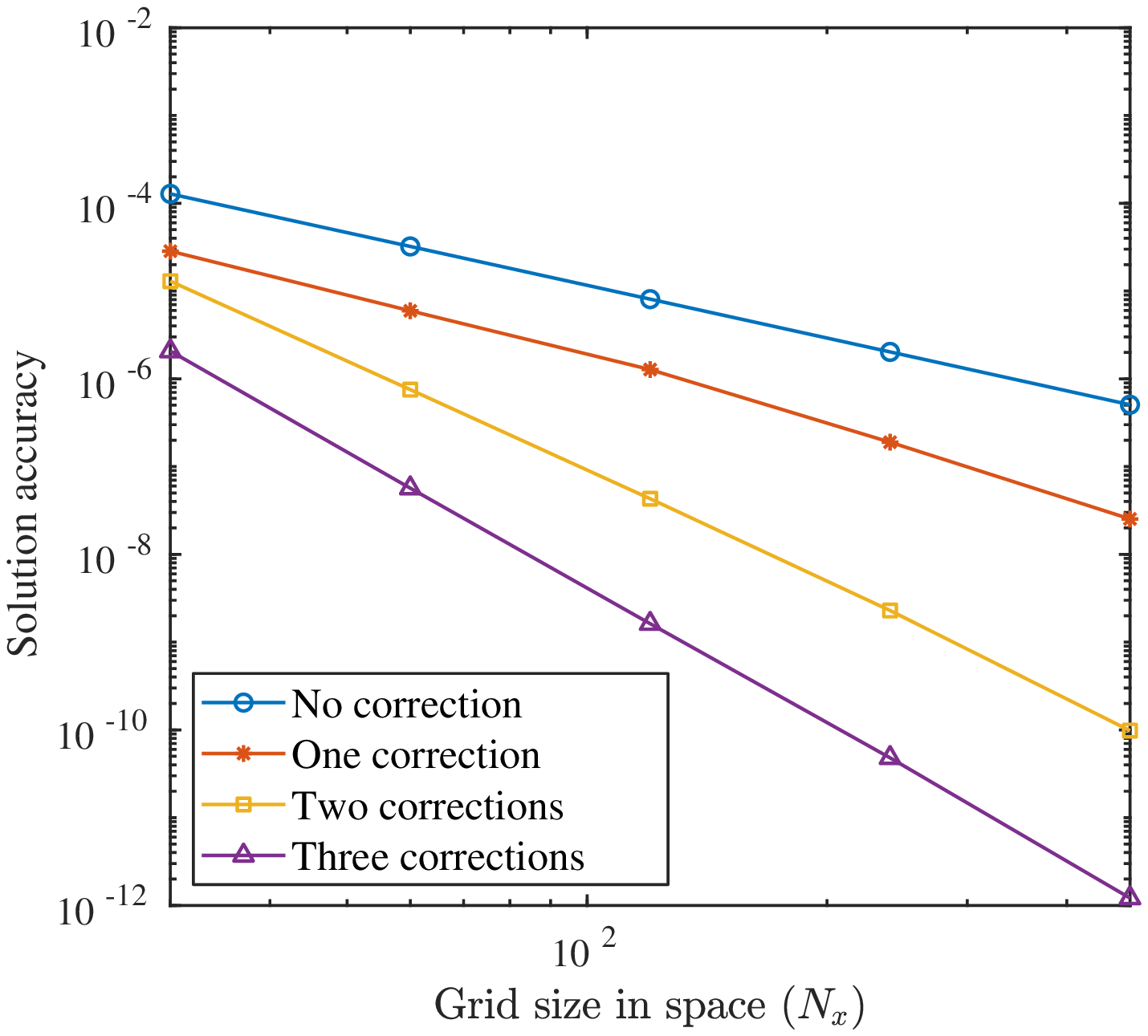}\\
    (a)  & (b)
\end{tabular}
    \caption{Log-log plot of solution error at point $x = 0.2$ versus computational complexity (a), and grid size in space (b), using results of Table \ref{tbl:ex1_bvp} for solving the penalty equation (\ref{eqn:bvp_ex_penalty}) of a one-dimensional free boundary obstacle problem. The computational complexity is represented by the grid size times the total number of penalty iterations.}
    \label{fig:bvp_complexity}
\end{figure}

We use this example to demonstrate that the four solve phases with deferred corrections in Algorithm \ref{alg:fb} improve the convergence rates as expected. In Table \ref{tbl:ex1_bvp}, we can see that, away from the free boundary $x = 0$, the convergence orders at the first, second, third and fourth solves are 2, 3, 4 and 5, respectively. The free boundary approximation also follows the same successive increase of convergence order, as shown in Table \ref{tbl:ex1_bvp_fb_approx}. To demonstrate the computational efficiency of our algorithm, we have plotted the log-log graph of the solution errors versus the computational complexity, represented by grid size in space multiplied by the total number of penalty iterations, as shown in Figure \ref{fig:bvp_complexity}. Note that in the first solve phase with no corrections, we only use a rough initial guess of a constant function $f = 1$, even though a better initial guess could be chosen. As a result, the first solve phase requires several penalty iterations to converge and takes up the major computational cost of the algorithm. In the second to fourth solve phases, only a single iteration is required for each solve phase, because the solution of the previous solve phase provides a good initial guess for Newton's method.

\begin{remark1}
Note that in both Tables \ref{tbl:ex1_bvp} and \ref{tbl:ex1_bvp_fb_approx}, we see that fifth-order convergence is obtained after applying three corrections. The reason for this is that, even though we are using a fourth-order method, the error from the nonsmoothness of the solution at the free boundary has been reduced to $\mathcal O(h^5)$. Since this error is large compared to the $\mathcal O(h^4)$ error arising from the fourth-order finite difference scheme, the convergence appears to be fifth-order.
\end{remark1}

\begin{table}[h]
\centering
\begin{tabular}{c|c|c|c|c|c|c|c|c}
\midrule\midrule
\multirow{3}{*}{N} & \multicolumn{8}{c}{$x = 0.2$}\\ \cline{2-9}
& \multicolumn{4}{c|}{1st solve (no correction)}& \multicolumn{4}{c}{2nd solve (one correction)}\\\cline{2-9}
 & \text{niters} & \text{value}  & \text{error} & \text{conv} & \text{niters} & \text{value}  & \text{error} & \text{conv}\\
\hline
30 & 7 & 0.221530668 & 1.28e-04 & - & 8 & 0.221431233 & 2.85e-05 & - \\
60 & 12 & 0.221434987 & 3.22e-05  & 1.99 & 13 & 0.221396803 & 5.96e-06 & 2.26 \\
120 & 23 & 0.221410846 & 8.09e-06 & 1.99 & 24 & 0.221401484 & 1.27e-06 & 2.22\\
240 & 44 & 0.221404784 & 2.03e-06 & 2.00 & 45 & 0.221402568 & 1.90e-07 & 2.75\\
480 & 86 & 0.221403265 & 5.07e-07 & 2.00 & 87 & 0.221402733 & 2.53e-08 & 2.91\\
\hline
\multirow{2}{*}{N} & \multicolumn{4}{c|}{3rd solve (two corrections)}& \multicolumn{4}{c}{4th solve (three corrections)}\\\cline{2-9}
& \text{niters} & \text{value}  & \text{error} & \text{conv} & \text{niters} & \text{value}  & \text{error} & \text{conv} \\
\hline
30 & 9 & 0.221415708 & 1.29e-05 & - & 10 & 0.221400689 & 2.07e-06 & -\\
60 & 14 & 0.221403511 & 7.53e-07 & 4.10 & 15 & 0.221402701 & 5.68e-08 & 5.19\\
120 & 25 & 0.221402801 & 4.30e-08 & 4.13 & 26 & 0.221402757 & 1.62e-09 & 5.13\\
240 & 46 & 0.221402760 & 2.29e-09 & 4.23 & 47 & 0.221402758 & 4.77e-11 & 5.09\\
480 & 88 & 0.221402758 & 9.83e-11 & 4.54 & 89 & 0.221402758 & 1.06e-12 & 5.50\\
\midrule\midrule
\end{tabular}
\caption{Convergence results of solutions at point $x=0.2$, of each solve phase in Algorithm \ref{alg:fb} for solving the penalty equation (\ref{eqn:bvp_ex_penalty}) of a one-dimensional free boundary obstacle problem with free boundary at $x = 0$. Uniform grid spacing is used. Note that ``niters" for the second to fourth solve includes the total number of iterations from all previous solve phases.}\label{tbl:ex1_bvp}
\end{table}

\begin{table}[h]
\centering
\begin{tabular}{c|c|c|c|c|c|c|c|c}
\midrule\midrule
\multirow{2}{*}{N} & \multicolumn{2}{c|}{1st solve}& \multicolumn{2}{c|}{2nd solve} & \multicolumn{2}{c|}{3rd solve}& \multicolumn{2}{c}{4th solve}\\\cline{2-9}
 & \text{error} & \text{conv}  & \text{error} & \text{conv} & \text{error} & \text{conv}  & \text{error} & \text{conv}\\
\hline
30 & 1.82e-02 & - & 2.20e-03 & - & 2.20e-04 & - & 1.87e-05 & -\\
60 & 4.69e-03 & 1.96 & 2.51e-04 & 3.13 & 1.15e-05 & 4.26 & 4.12e-07 & 5.50\\
120 & 1.19e-03 & 1.98 & 2.84e-05 & 3.15 & 5.85e-07 & 4.30 & 6.55e-09 & 5.97\\
240 & 2.97e-04 & 2.00 & 2.93e-06 & 3.28 & 2.58e-08 & 4.50 & 8.38e-11 & 6.29\\
480 & 7.33e-05 & 2.02 & 2.17e-07 & 3.75 & 1.02e-09 & 4.67 & 1.03e-11 & 3.02\\
\midrule\midrule
\end{tabular}
\caption{Convergence results of the free boundary approximation of each solve phase in Algorithm \ref{alg:fb} for solving the penalty equation (\ref{eqn:bvp_ex_penalty}) of a one-dimensional free boundary obstacle problem with the exact free boundary at $x = 0$. Uniform grid spacing is used.}\label{tbl:ex1_bvp_fb_approx}
\end{table}

\subsection{A simple test moving boundary problem}\label{sec:simple-ivp-numerical-results}
In this second example, we introduce the time variable and consider a time-dependent free boundary problem. Consider the LCP
\begin{equation}\label{eqn:ivp_freeb_ex}
\begin{aligned}
f_t-\frac{1}{2\sqrt{t}}f'' + \frac{1}{2\sqrt{t}}&\geq 0,\\
f - f^* &\geq 0,\\
\left(f_t-\frac{1}{2\sqrt{t}}f'' + \frac{1}{2\sqrt{t}} = 0\right) \vee (f-f^*&=0),
\end{aligned}
\end{equation}
on the domain $(t,x)\in[0,0.5]\times[-2,2]$, where $f^*(t,x) = x$. The solution satisfies the Dirichlet boundary conditions
\begin{equation*}
    f(t,-2) = -2,\; f(t,2) = e^{2+\sqrt{t}}-\sqrt{t}-1,
\end{equation*}
and the initial condition
\begin{equation*}
    f(0,x) = \left\{
    \begin{array}{lr}
    e^{x} - 1, & 0 \leq x \leq 2,\\
    x, & -2 \leq x < 0.
    \end{array}\right.
\end{equation*}
The exact solution to (\ref{eqn:ivp_freeb_ex}) is
\begin{equation*}
    f(t,x) = \left\{
    \begin{array}{lr}
    e^{x+\sqrt{t}} - \sqrt{t} - 1, & x_f(t) \leq x \leq 2,\\
    x, & -2 \leq x < x_f(t),
    \end{array}\right.
\end{equation*}
where $x_f(t)$ is the moving free boundary
\begin{equation*}
    x_f(t) = -\sqrt{t}.
\end{equation*}
The value matching and smooth pasting conditions at the free boundary $x = x_f(t)$ follow naturally. Again, we see that $f(\cdot,x)\in C^1\backslash C^2$ on $[-2,2]$, but it is smooth on $[-2,x_f(t)]$ and $[x_f(t), 2]$, separately. To apply Algorithm \ref{alg:mb}, we write (\ref{eqn:ivp_freeb_ex}) in penalty form
\begin{equation}\label{eqn:ivp_freeb_ex_penalty}
    f_t-\frac{1}{2\sqrt{t}}f'' + \frac{1}{2\sqrt{t}} + \rho\max(f^*-f,0) = 0,
\end{equation}
where $\rho$ is a penalty constant, taken to be $\rho = 1\times 10^8$ in the numerical experiments.

Since the free boundary $x_f(t) = -\sqrt{t}$, its location changes rapidly near time $t = 0$. This will cause a problem in the BDF4 time-stepping scheme because many grid points will cross the free boundary in the initial time steps (see Section \ref{sec:time-discretization}). Hence, BDF4 degenerates to only first-order convergence due to piecewise smoothness in the solution across the free boundary. To avoid this situation, we perform a time-variable transformation $t = \tau^2$ so that the free boundary changes more slowly, and fewer points will cross the free boundary in the initial time steps. Although this does not completely solve the problem, it is accurate enough for the algorithm to achieve high-order convergence, as shown in the numerical results. To start BDF4, we use the exact solutions for the first three time steps. For this problem, we simply use a uniform grid in space.

\begin{figure}[h]
    \centering
\begin{tabular}{c c}
    \includegraphics[width=0.48\textwidth]{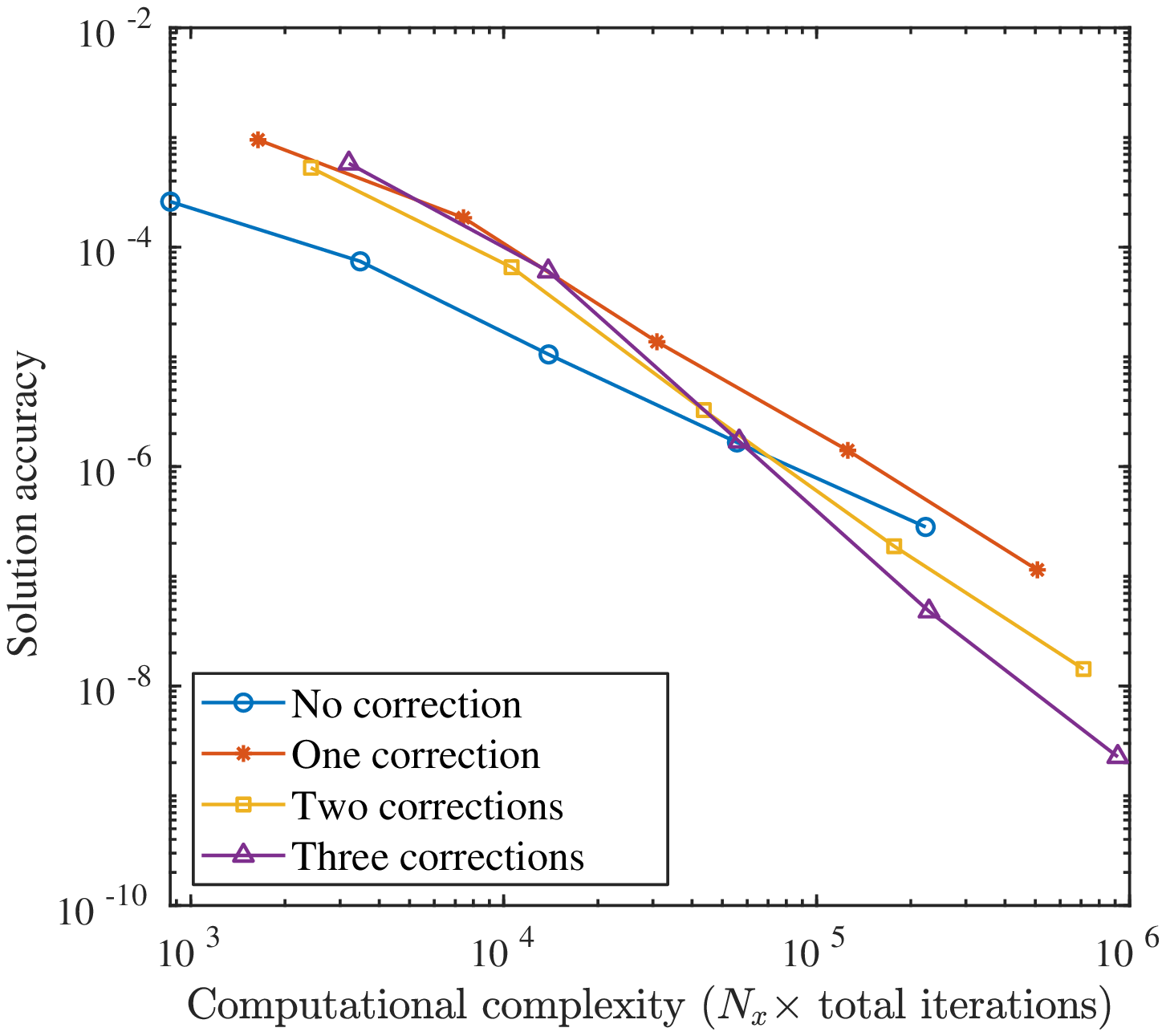} &  \includegraphics[width=0.48\textwidth]{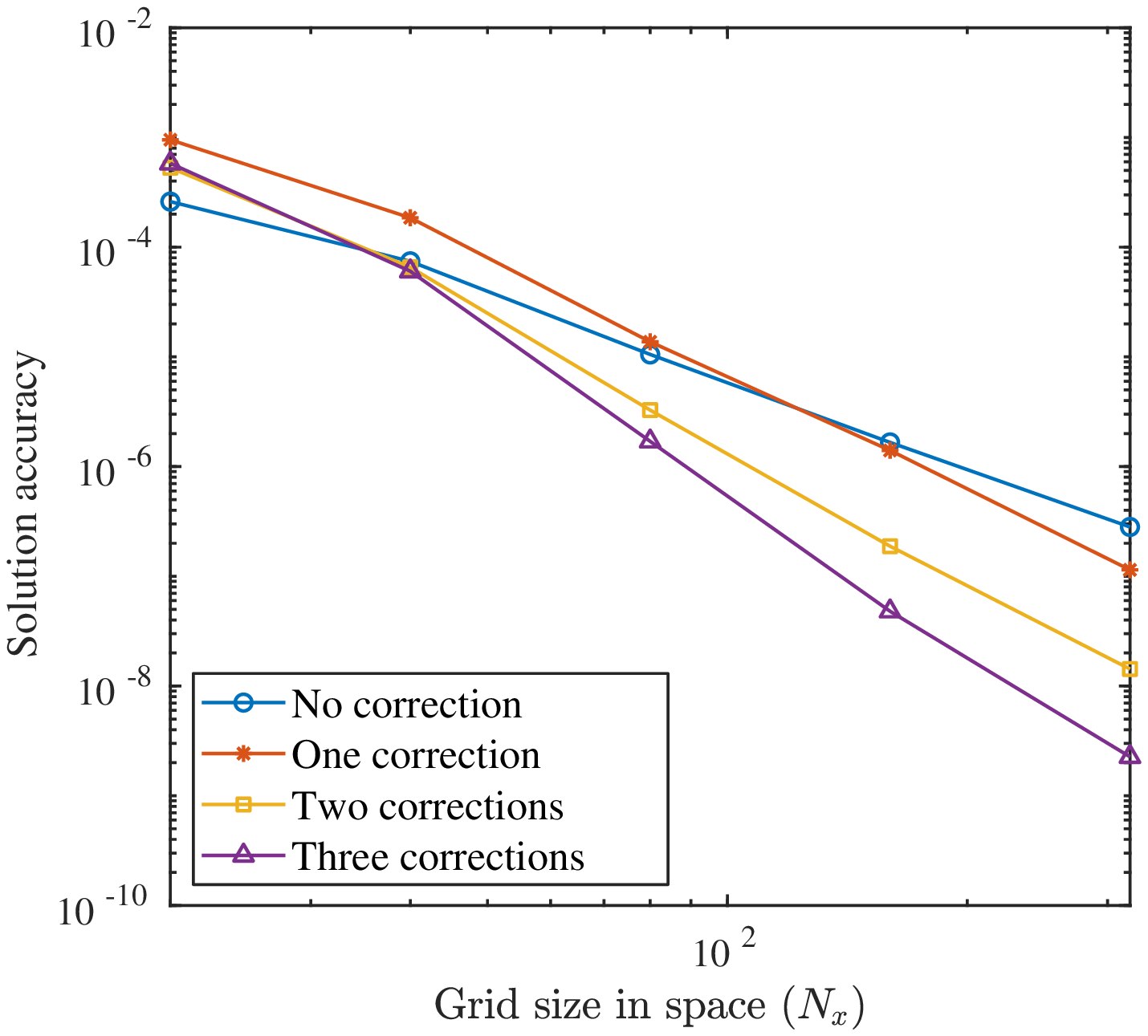}\\
    (a)  & (b)
\end{tabular}
    \caption{Log-log plot of solution errors at point $x = 0$ versus computational complexity (a), and grid size in space (b), using results of Table \ref{tbl:ex2_freeb_1} for solving the penalty equation (\ref{eqn:ivp_freeb_ex_penalty}) of a moving boundary problem with the exact moving boundary $x_f(t) = -\sqrt{t}$. The computational complexity is represented by the grid size times the total number of penalty iterations.}
    \label{fig:smooth_fb_complexity}
\end{figure}

In Table \ref{tbl:ex2_freeb_1}, we record the convergence results at $x = -0.37$ and at $x = 0$. The point $x = -0.37$ is slightly to the right of the first grid point right of the final-time free boundary location on the coarsest grid $N_x = 20$. The point $x = 0$ is the initial free boundary location. We see that the solutions at both points gain the expected order of convergence after each correction. To demonstrate the computational efficiency of our algorithm, we plot the log-log graph of the solution errors versus the computational complexity represented by the grid size in space multiplied by the total number of penalty iterations, as shown in Figure \ref{fig:smooth_fb_complexity}.

\begin{remark1}
Note that the solutions after solving with corrections have larger errors on the coarsest grid $N_x = 20$. This is due to large extrapolation errors of free boundary and derivatives approximations when the space step size near the free boundary is large, which occurs on a uniform coarse grid. This can be avoided by applying grid stretching around the free boundary.
\end{remark1}

\begin{table}
\centering
\begin{tabular}{c|c|c|c|c|c|c|c|c}
\midrule\midrule
\multirow{3}{*}{$(N_x, N_t)$} & \multicolumn{8}{c}{$x = -0.37$}\\ \cline{2-9}
& \multicolumn{4}{c|}{1st solve (no correction)}& \multicolumn{4}{c}{2nd solve (one correction)}\\\cline{2-9}
 & \text{niters} & \text{value}  & \text{error} & \text{conv} & \text{niters} & \text{value}  & \text{error} & \text{conv}\\
\hline
(20,40) & 43 & -0.307227 & 1.04e-03 & - & 82 & -0.304930 & 1.26e-03 & -\\
(40,80) & 87 & -0.306241 & 8.16e-06 & 7.00 & 186 & -0.306075 & 1.47e-04 & 3.10\\
(80,160) & 174 & -0.306227 & 9.13e-06 & -0.16 & 386 & -0.306205 & 1.33e-05 & 3.46\\
(160,320) & 348 & -0.306220 & 1.94e-06 & 2.24 & 786 & -0.306216 & 1.67e-06 & 3.00\\
(320,640) & 697 & -0.306218 & 3.19e-07 & 2.60 & 1586 & -0.306218 & 1.44e-07 & 3.53\\
\hline
\multirow{2}{*}{$(N_x, N_t)$} & \multicolumn{4}{c|}{3rd solve (two corrections)}& \multicolumn{4}{c}{4th solve (three corrections)}\\\cline{2-9}
& \text{niters} & \text{value}  & \text{error} & \text{conv} & \text{niters} & \text{value}  & \text{error} & \text{conv} \\
\hline
(20,40) & 121 & -0.306334 & 1.61e-04 & - & 160 & -0.305883 & 3.03e-04 & -\\
(40,80) & 265 & -0.306147 & 7.44e-05 & 1.11 & 347 & -0.306163 & 5.90e-05 & 2.36\\
(80,160) & 545 & -0.306215 & 3.35e-06 & 4.47 & 707 & -0.306216 & 1.77e-06 & 5.06\\
(160,320) & 1106 & -0.306218 & 2.20e-07 & 3.93 & 1429 & -0.306218 & 5.78e-08 & 4.94\\
(320,640) & 2227 & -0.306218 & 1.12e-08 & 4.29 & 2866 & -0.306218 & 2.79e-09 & 4.37\\
\hline
\multirow{3}{*}{$(N_x, N_t)$} & \multicolumn{8}{c}{$x = 0$}\\ \cline{2-9}
& \multicolumn{4}{c|}{1st solve (no correction)}& \multicolumn{4}{c}{2nd solve (one correction)}\\\cline{2-9}
 & \text{niters} & \text{value}  & \text{error} & \text{conv} & \text{niters} & \text{value}  & \text{error} & \text{conv}\\
\hline
(20,40) & 43 & 0.320748 & 2.60e-04 & - & 82 & 0.321960 & 9.52e-04 & -\\
(40,80) & 87 & 0.320934 & 7.43e-05 & 1.81 & 186 & 0.321194 & 1.85e-04 & 2.36\\
(80,160) & 174 & 0.320998 & 1.05e-05 & 2.82 & 386 & 0.321022 & 1.37e-05 & 3.76\\
(160,320) & 348 & 0.321007 & 1.66e-06 & 2.66 & 786 & 0.321010 & 1.40e-06 & 3.29\\
(320,640) & 697 & 0.321008 & 2.82e-07 & 2.56 & 1586 & 0.321008 & 1.14e-07 & 3.62\\
\hline
\multirow{2}{*}{$(N_x, N_t)$} & \multicolumn{4}{c|}{3rd solve (two corrections)}& \multicolumn{4}{c}{4th solve (three corrections)}\\\cline{2-9}
& \text{niters} & \text{value}  & \text{error} & \text{conv} & \text{niters} & \text{value}  & \text{error} & \text{conv} \\
\hline
(20,40) & 121 & 0.321537 & 5.29e-04 & - & 160 & 0.321588 & 5.80e-04 & -\\
(40,80) & 265 & 0.321074 & 6.56e-05 & 3.01 & 347 & 0.321069 & 6.04e-05 & 3.26\\
(80,160) & 545 & 0.321011 & 3.27e-06 & 4.33 & 707 & 0.321010 & 1.70e-06 & 5.15\\
(160,320) & 1106 & 0.321008 & 1.88e-07 & 4.12 & 1429 & 0.321008 & 4.80e-08 & 5.15\\
(320,640) & 2227 & 0.321008 & 1.43e-08 & 3.72 & 2866 & 0.321008 & 2.27e-09 & 4.41\\
\midrule\midrule
\end{tabular}
\caption{Convergence results of solutions at points $x = -0.37$ and $x=0$ of each solve phase in Algorithm \ref{alg:mb} for solving the penalty equation (\ref{eqn:ivp_freeb_ex_penalty}) of a moving boundary problem with the exact moving boundary $x_f(t) = -\sqrt{t}$. Note that ``niters" for the second to fourth solve includes the total number of iterations from all previous solve phases. }\label{tbl:ex2_freeb_1}
\end{table}

\subsection{American option pricing}\label{sec:american-numerical-results}
Finally, we use our algorithm to solve the American put option pricing problem. We repeat the penalty equation
\begin{equation*}
    \partial_t V = \mathcal L_{BS}V + \rho \max\{V^* - V, 0\},
\end{equation*}
for convenience of discussion, where $V^*(S) = \max\{K-S, 0\}$ is the payoff of the American put option struck at $K$.
The initial condition is
\begin{equation*}
    V(0,S) = V^*(S).
\end{equation*}
We truncate the right end of the domain at $S = S_{\max}$ and use Dirichlet boundary conditions $V(t,0) = K$ and $V(t,S_{\max}) = 0$.
To avoid complications due to the payoff singularity of the first derivative at the strike price of the American put options, we compute the difference between an American option and a European option, as in \cite{zhu2004derivative}, where this is referred to as the singularity-separating method.
A European put option value $V^{E}$ with the same volatility $\sigma$, bank interest $r$, dividend $q$, and strike price $K$ satisfies the Black-Scholes equation $\partial_t V^E = \mathcal L_{BS}V^E$,
and the initial condition $V^E(0,S) = \max\{K-S,0\}$, and has a known explicit formula.
The difference of the solutions $V^{\text{diff}}=V-V^E$ has a zero initial condition.
Therefore, instead of solving for the original American option price, we solve for $V^{\text{diff}}$, which satisfies the equation
\begin{equation*}
    \partial_t V^{\text{diff}} = \mathcal L_{BS}V^{\text{diff}} + \rho \max\{(V^*-V^E) - V^{\text{diff}}, 0\},
\end{equation*}
with a zero initial condition, and then add the European option price back to obtain the final American option price. The penalty constant $\rho$ is chosen to be $\rho = 1\times10^8$ in the numerical experiments.

In this problem, we do not have the exact solution. Since BDF4 requires solutions from the previous four time steps to proceed, we need to be careful when starting BDF4. 
For this problem, we compute the first time step using the classical fourth-order Runge-Kutta method \cite{butcher1996history}. For the second time step, we use a three-level fourth-order method
\begin{equation*}
    \left(\mathbf I - \frac{k}{3}\mathbf L\right)\tilde{\mathbf V}^{n+2} = \tilde{\mathbf V}^n + \frac{k}{3}\mathbf L\left( 4\tilde{\mathbf V}^{n+1} + \tilde{\mathbf V}^n \right) + \frac{k}{3}(\mathbf b^{n+2} + 4\mathbf b^{n+1} + \mathbf b^n) + 2\mathbf q(\tilde{\mathbf V}^{n+2}),
\end{equation*}
see \cite{li2001compact}, and for the third time step, we simply apply BDF3. We observe that this starting scheme is sufficient for fourth-order convergence.

We test the algorithm on two example problems with different volatilities, $\sigma = 0.2$ and $\sigma = 0.8$ as the examples in \cite{forsyth2002quadratic}.
Both examples have the other parameters the same: zero dividend payment,
interest rate $r = 0.1$, strike price $K = 100$ and expiration time $T = 0.25$.
We truncate the infinite domain at $S_{\max} = 10K = 1000$ for the problem with smaller volatility $\sigma = 0.2$,
and at $S_{\max}=13K = 1300$ for the larger volatility $\sigma = 0.8$.
As it turns out, a larger $S_{\max}$ is not only necessary for the accuracy of solution with a larger volatility,
it is also important for observing the convergence results of our algorithm.
To implement the algorithm, we use the time variable transformation $t = \tau^2$ for both examples.
We propose the stretching function
\begin{equation*}
    \xi(S) = \left(S - \frac{\sqrt{\pi}}{2}\frac{1-\beta}{\beta}\alpha\text{erfc}\left(  \frac{S-K}{\alpha}\right)\right)C_1 + C_2,
\end{equation*}
to stretch the space grid, where
\begin{align*}
    &C_1 = 1/\left[ (S_{\max}-S_{\min}) - \frac{\sqrt{\pi}}{2}\frac{1-\beta}{\beta}\alpha\left( \text{erfc}\left( \frac{S_{\max}-K}{\alpha} \right)  - \text{erfc}\left( \frac{S_{\min}-K}{\alpha} \right) \right) \right],\\
    &C_2 = \left[ \frac{\sqrt{\pi}}{2}\frac{1-\beta}{\beta}\alpha\text{erfc}\left( \frac{S_{\min}-K}{\alpha} \right) - S_{\min} \right]C_1.
\end{align*}
where $\alpha$ and $\beta$ are parameters controlling the density of the stretching, and chosen to be $\alpha = 125/6$ and $\beta = 1/20$ for $\sigma = 0.2$, and $\alpha = 65, \beta = 1/8$ for $\sigma = 0.8$, in the numerical experiments. This function adds additional grid points with density $1/\beta$, on a region of width $6\alpha$ centered around the point $K$, while maintaining the density of 1 elsewhere. Our choice of parameters $\alpha$ can be understood by the fact that: for $\sigma = 0.2$, the range of the optimal exercise boundary movement is approximately 10, starting from $S = 100$ at $t = 0$, to $S = 89.7$ at expiry $t = 0.25$. Therefore, the moving free boundary is within $\frac{10}{3\alpha}\approx\frac{1}{6}$ of the length of the stretched region, away from the stretching center, during the whole time period. For $\sigma = 0.8$, the range of the optimal exercise boundary movement is around 48, starting from $S = 100$ at $t = 0$, to $S = 51.8$ at expiry $t = 0.25$. Therefore, the moving free boundary is within $\frac{51.8}{3\alpha} \approx \frac{1}{4}$ of the length of the stretched region, away from the stretching center, during the whole time period.

In addition, since the solution of the American option price has singular derivative at the strike at expiry, meaning that the solution is not smooth, we do not apply the correction scheme in Algorithm \ref{alg:mb} for the first several time steps. In the numerical experiments, the number of time steps skipped is chosen to be $t_{skip} = 12$ for both examples. Since we apply time and space stretching near the strike at expiry, the errors from the skipped corrections in the first few time steps are sufficiently small so as to not affect the high-order convergence.

In Table \ref{tbl:ex4_american_1}, we can see clear convergence rate improvements and error reduction with corrections at each solve phase, with the smaller volatility exhibiting even faster error reduction. The second solve phase, corresponding to one correction, with first correction only slightly changes the error. The third- and fourth-solve phases, corresponding to two and three corrections, respectively, reduce the error significantly. The final results after the fourth solve phase exhibit a reduction of error by nearly 100 times compared to the no-correction phase. 

\begin{figure}[h]
    \centering
\begin{tabular}{c c}
    \includegraphics[width=0.48\textwidth]{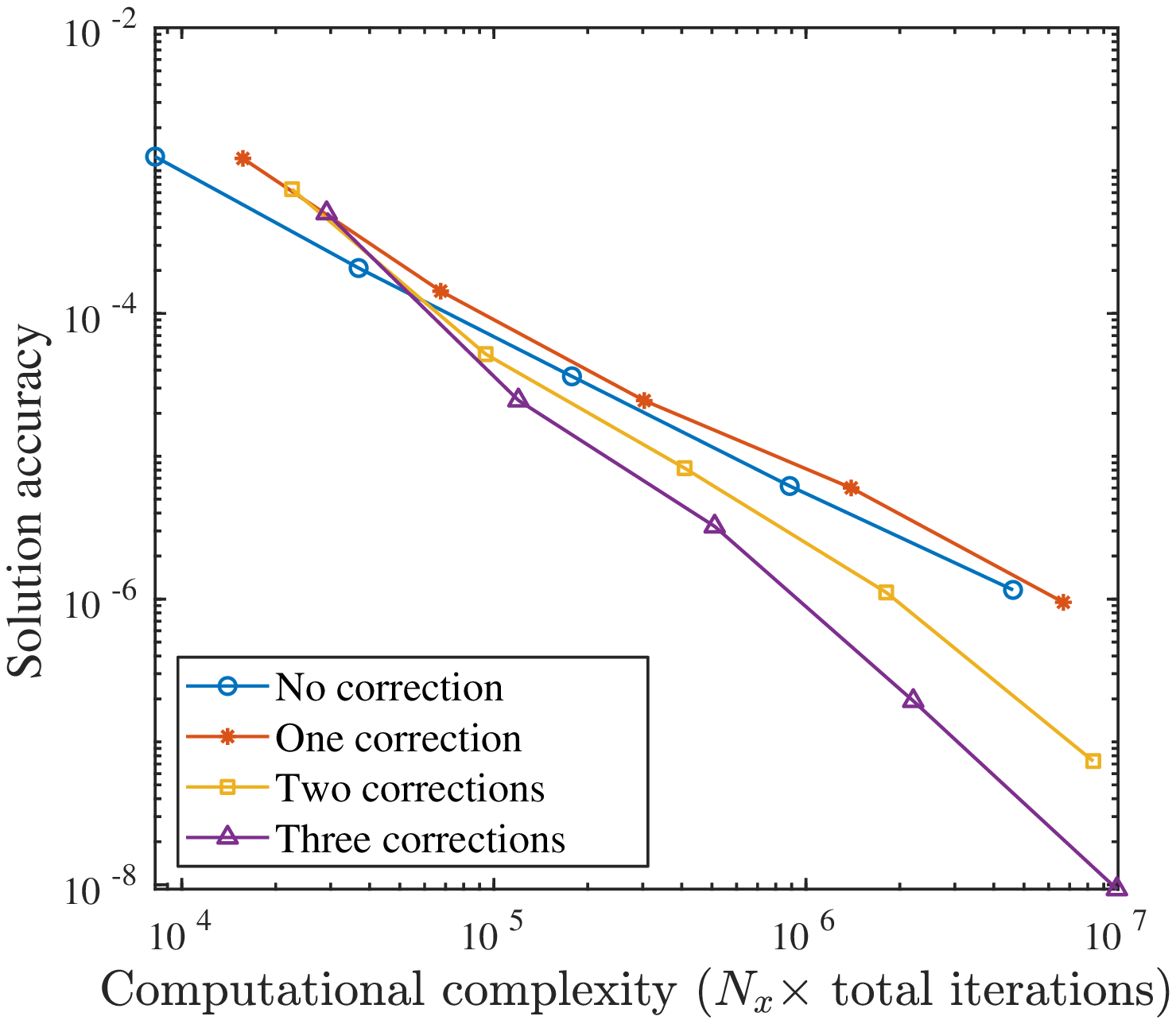} &  \includegraphics[width=0.47\textwidth]{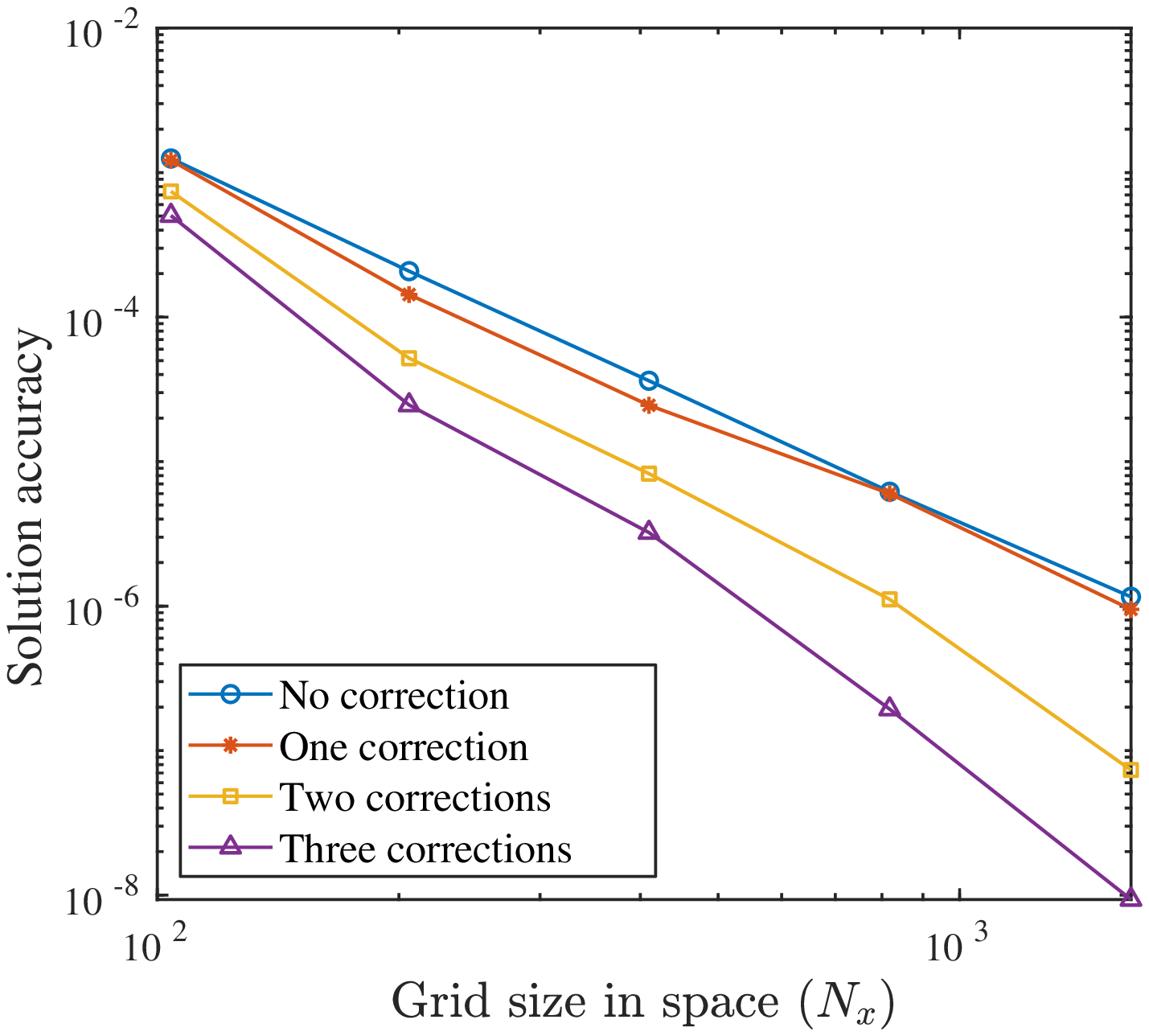}\\
    (a)  & (b)
\end{tabular}
    \caption{Log-log plot of solution changes at the strike point $K$ at the final time $T$ versus computational complexity (a), and grid size in space (b), using results of solving American option prices in Table \ref{tbl:ex4_american_1} for $\sigma = 0.2$. The computational complexity is represented by the grid size in space times the total number of penalty iterations.}
    \label{fig:am_sig02_complexity}
\end{figure}

To demonstrate the computational efficiency of our algorithm, we have shown in Figure \ref{fig:am_sig02_complexity} the solution accuracy versus the computational complexity represented by the grid size in space multiplied with the total number of penalty iterations. We can see that the three-correction algorithm is slightly more expensive if high accuracy is not the goal.
However, when a high accuracy solution is desired, the three-correction algorithm is more efficient.

\begin{remark1}
For this example, we see that the second solve does not improve the convergence much. This is because we have applied enough stretching to reduce the leading-order error term due to the second-derivative jump in the first solve, even without correction. Moreover, we observe that our algorithm for $\sigma = 0.2$ performs better than for $\sigma = 0.8$. For a larger volatility $\sigma = 0.8$, the optimal exercise boundary moves more quickly and ranges over a larger part of the domain within the same time span. Since we only apply stretching in space around the initial free boundary, when the free boundary moves farther away from the stretching center, the extrapolation error when approximating the free boundary and derivatives becomes larger due to larger grid spacing. This explains the increase in error from applying the first correction. One way to avoid large extrapolation errors is to implement a time-dependent grid stretching that follows the free boundary movement, see e.g. \cite{oosterlee2005accurate} where a predictor-corrector scheme is applied. We leave this to future work.
\end{remark1}

\begin{figure}
    \centering
\begin{tabular}{c c}
    \includegraphics[width=0.48\textwidth]{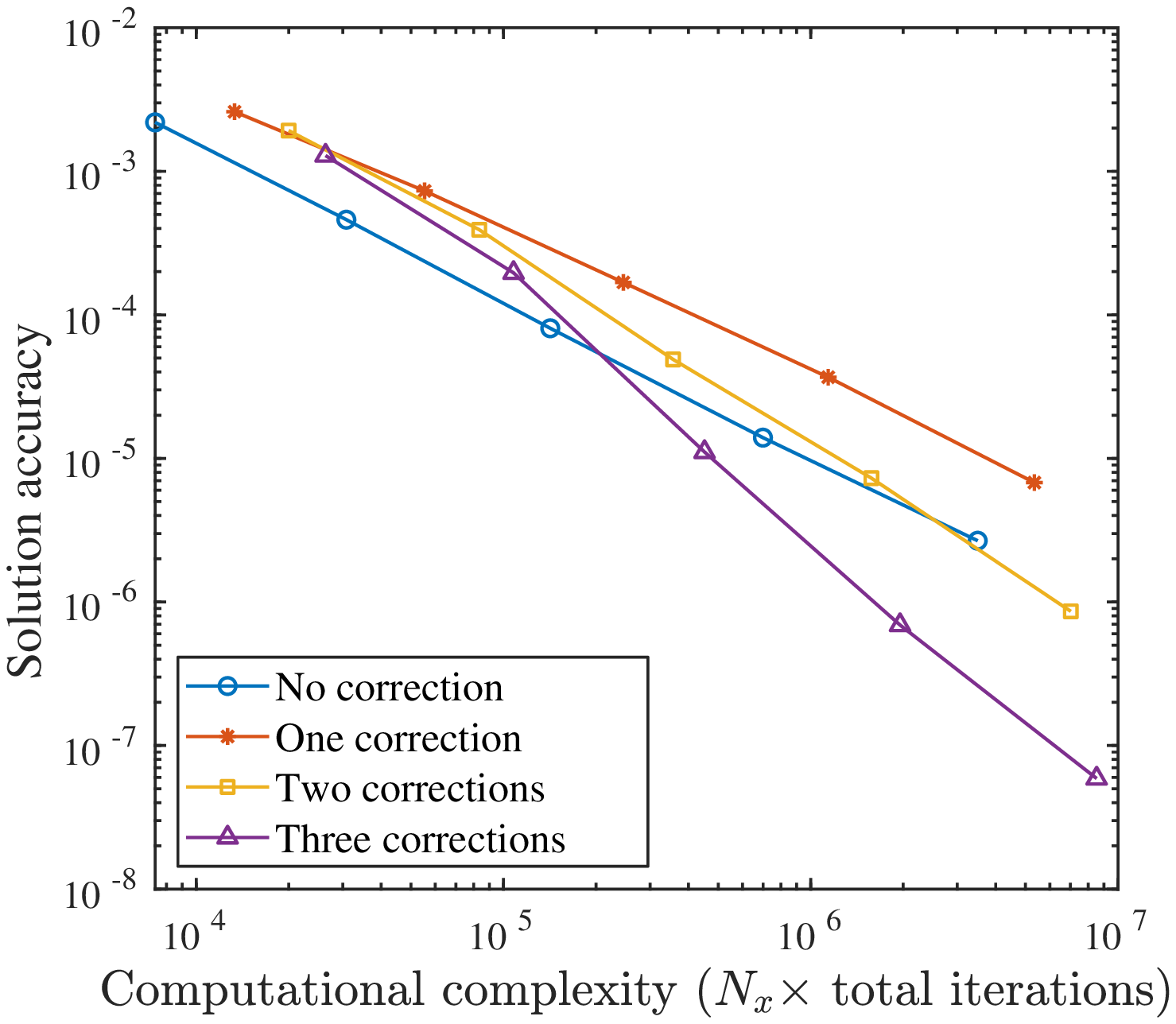} &  \includegraphics[width=0.47\textwidth]{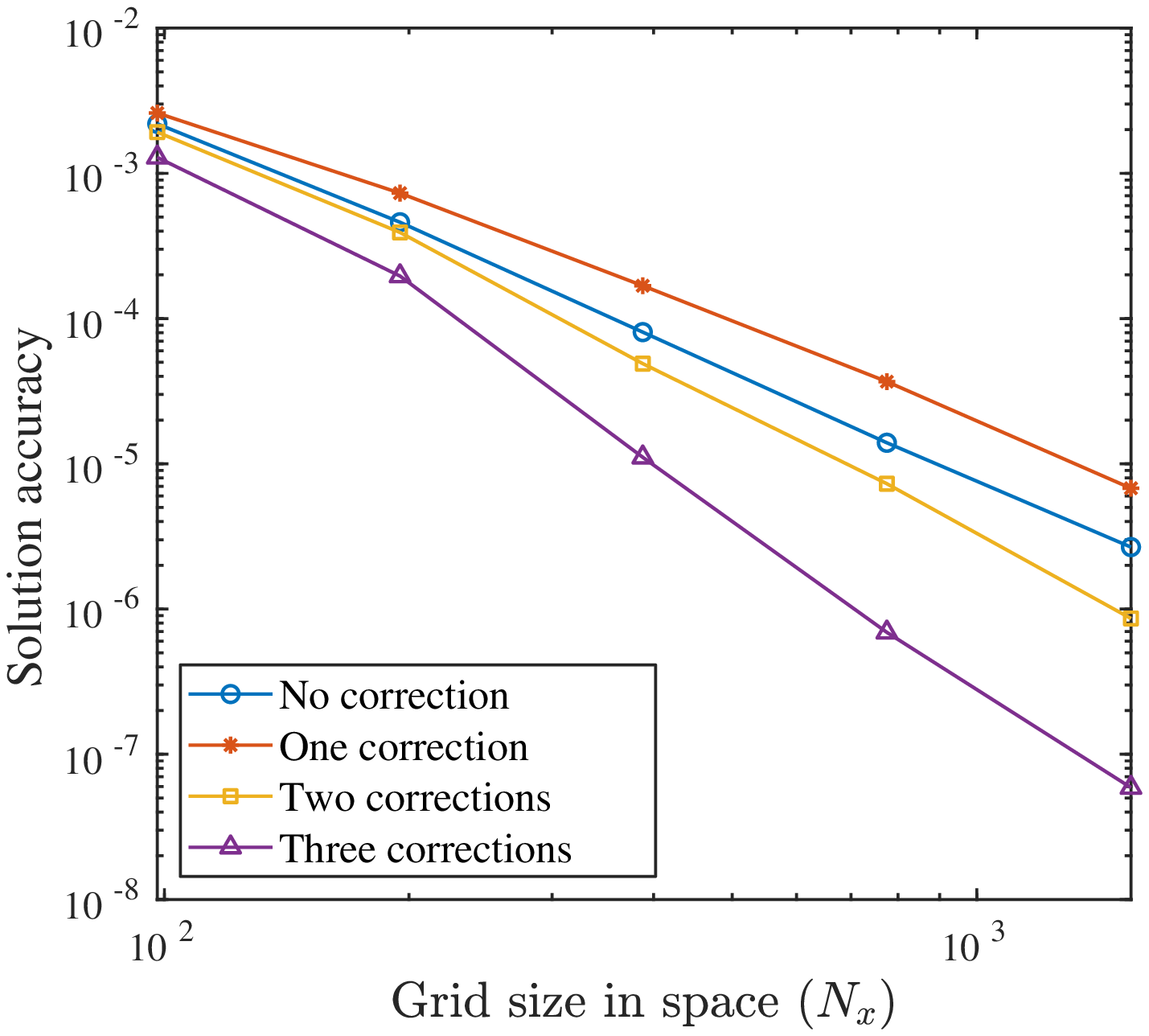}\\
    (a)  & (b)
\end{tabular}
    \caption{Log-log plot of solution changes at the strike point $K$ at the final time $T$ versus computational complexity (a), and grid size in space (b), using results of solving American option prices in Table \ref{tbl:ex4_american_1} for $\sigma = 0.8$. The computational complexity is represented by the grid size in space times the total number of penalty iterations.}
    \label{fig:am_sig08_complexity}
\end{figure}

\begin{table}
\centering
\begin{tabular}{c|c|c|c|c|c|c|c|c}
\midrule\midrule
\multirow{3}{*}{$(N_x, N_t)$} & \multicolumn{8}{c}{$\sigma = 0.2,\; T = 0.25$}\\ \cline{2-9}
& \multicolumn{4}{c|}{1st solve (no correction)}& \multicolumn{4}{c}{2nd solve (one correction)}\\\cline{2-9}
 & \text{niters} & \text{value}  & \text{error} & \text{conv} & \text{niters} & \text{value}  & \text{error} & \text{conv}\\
\hline
(53,30) & 37 & 3.068602382 & - & - & 71 & 3.068715191 & - & -\\
(104,60) & 79 & 3.069855016 & 1.25e-03 & - & 151 & 3.069931750 & 1.22e-03 & -\\
(206,120) & 179 & 3.070062874 & 2.08e-04 & 2.59 & 328 & 3.070075013 & 1.43e-04 & 3.09\\
(410,240) & 434 & 3.070099140 & 3.63e-05 & 2.52 & 740 & 3.070099544 & 2.45e-05 & 2.55\\
(818,480) & 1085 & 3.070105329 & 6.19e-06 & 2.55 & 1708 & 3.070105547 & 6.00e-06 & 2.03\\
1635,960) & 2820 & 3.070106489  & 1.16e-06 & 2.42 & 4088 & 3.070106496 & 9.49e-07 & 2.66\\
\hline
\multirow{2}{*}{$(N_x, N_t)$} & \multicolumn{4}{c|}{3rd solve (two corrections)}& \multicolumn{4}{c}{4th solve (three corrections)}\\\cline{2-9}
& \text{niters} & \text{value}  & \text{error} & \text{conv} & \text{niters} & \text{value}  & \text{error} & \text{conv} \\
\hline
(53,30) & 104 & 3.069304376 & - & - & 136 &  3.069574295 & - & -\\
(104,60) & 217 & 3.070045347 & 7.41e-04 & - & 280 & 3.070078659 & 5.04e-04 & -\\
(206,120) & 457 & 3.070097280 & 5.19e-05 & 3.83 & 582 & 3.070103304 & 2.46e-05 & 4.36\\
(410,240) & 997 & 3.070105534 & 8.25e-06 & 2.65 & 1244 & 3.070106531 & 3.23e-06 & 2.93\\
(818,480) & 2211 & 3.070106647 & 1.11e-06 & 2.89 & 2699 &  3.070106725 & 1.94e-07 & 4.06\\
(1635,960) & 5096 &  3.070106721 & 7.34e-08 & 3.92 & 6067 & 3.070106734 & 9.32e-09 & 4.38\\
\hline
\multirow{3}{*}{$(N_x, N_t)$} & \multicolumn{8}{c}{$\sigma = 0.8,\; T = 0.25$}\\ \cline{2-9}
& \multicolumn{4}{c|}{1st solve (no correction)}& \multicolumn{4}{c}{2nd solve (one correction)}\\\cline{2-9}
 & \text{niters} & \text{value}  & \text{error} & \text{conv} & \text{niters} & \text{value}  & \text{error} & \text{conv}\\
\hline
(50,30) & 34 & 14.676127404 & - & - & 64 & 14.675332474 & - & -\\
(98,60) & 75 & 14.678320134 & 2.19e-03 & - & 136 & 14.677934361 & 2.60e-03 & -\\
(195,120) & 158 & 14.678780547 & 4.60e-04 & 2.25 & 284 & 14.678665160 & 7.31e-04 & 1.83\\
(388,240) & 366 & 14.678861193 & 8.06e-05 & 2.51 & 633 & 14.678833475 & 1.68e-04 & 2.12\\
(775,480) & 902 & 14.678875150 & 1.40e-05 & 2.53 & 1471 &  14.678870259 & 3.68e-05 & 2.19\\
(1548,960) & 2258 & 14.678877820 & 2.67e-06 & 2.39 & 3453 & 14.678877045 & 6.79e-06 & 2.44\\
\hline
\multirow{2}{*}{$(N_x, N_t)$} & \multicolumn{4}{c|}{3rd solve (two corrections)}& \multicolumn{4}{c}{4th solve (three corrections)}\\\cline{2-9}
& \text{niters} & \text{value}  & \text{error} & \text{conv} & \text{niters} & \text{value}  & \text{error} & \text{conv} \\
\hline
(50,30) & 94 & 14.676507458 & - & - & 126 & 14.677382768 & - & -\\
(98,60) & 204 &14.678429412 & 1.92e-03 & - & 269 & 14.678670286 & 1.29e-03 & -\\
(195,120) & 428 & 14.678821171 & 3.92e-04 &  2.29 & 553 &  14.678866522 & 1.96e-04 & 2.71\\
(388,240) & 919 & 14.678870105 & 4.89e-05 & 3.00 & 1163 &    14.678877609 & 1.11e-05 & 4.15\\
(775,480) & 2033 & 14.678877396 & 7.29e-06 & 2.75 & 2520 & 14.678878300 & 6.91e-07 & 4.00\\
(1548,960) & 4536 & 14.678878257 & 8.61e-07 & 3.08 & 5505 & 14.678878359 & 5.89e-08 & 3.55\\
\midrule\midrule
\end{tabular}
\caption{Convergence results of an American put option at $S = 100$, $T = 0.25$ with $K = 100,\; r = 0.1,\; q = 0$ for $\sigma = 0.2$ and $\sigma = 0.8$. Note that ``niters" for the second to fourth solve includes the total number of iterations from all previous solve phases.}\label{tbl:ex4_american_1}
\end{table}

\section{Conclusions}\label{sec:conclusion}
In this paper, we presented an analysis of the error when using fourth-order finite differences and BDF4 time-stepping scheme for solving free and moving boundary problems. Based on the analysis, we presented
a high-order deferred correction algorithm for solving these problems. Our algorithm utilizes the penalty method and assumes no prior knowledge of the exact free boundary location and derivative jumps at the free boundary.
Our method does not modify the finite difference stencils and the arising matrix, but applies the corrections to the right-hand side. The penalty iteration converges in a few iterations.
From the analysis of the error behaviors when solving free boundary problems, we showed that our deferred correction algorithm can successively increase the solution order of convergence from $\mathcal O(h^2)$ to $\mathcal O(h^3)$, and from $\mathcal O(h^3)$ to $\mathcal O(h^4)$ after applying each successive correction. Our numerical results validate the theoretical analysis. On simple test problems when the solution is not singular, results show that the behavior of our algorithm matches exactly with our theory. When solving the more challenging American put option problem, our algorithm also performs well.

\subsection{Generalizations and future work}
We only considered problems with one space dimension in this work, however, the deferred correction idea can be generalized to two space dimensions. One possible extension is to the elliptic obstacle problem in two dimensions, which is still an active area of research.
The major difference in the algorithm when applied to two-dimensional problems is how the extrapolation scheme is designed. Two dimensional extrapolation has already been extensively studied in the literature (see, for example, \cite{gibou2005fourth,linnick2005high,wiegmann2000explicit}). To reduce the extrapolation error on a uniform grid, we can also refine the grid around the free boundary, the location of which can be approximated on a coarser grid.

Since grid stretching in space is useful for reducing the extrapolation errors, another possible extension of our work is to apply a time-dependent grid stretching scheme when solving moving boundary problems. As noted in Section \ref{sec:results}, our algorithm applies grid stretching only around the initial free boundary, which can cause more extrapolation errors in later time steps when the free boundary moves out of the stretched area. Hence, it does not achieve ideal performance, as seen in the application to the American option pricing problem when $\sigma = 0.8$. By adapting the grid stretching so that it follows the moving boundary, the extrapolation error can be reduced at all time steps. One possible way of doing this is to use the predictor-corrector idea, as discussed in \cite{oosterlee2005accurate}. This scheme precomputes an approximate solution, which then gives an approximate moving boundary, on a coarse grid. When we solve the problem on a finer grid at each time step, we can apply grid stretching around the precomputed approximate free boundary. By doing this, we expect our algorithm to achieve even better results.

\appendix

\bibliographystyle{siam}
\bibliography{HOFB}

\end{document}